\magnification 1200
\baselineskip 14pt
\parskip=3pt plus1pt minus.5pt
\input amssym.def

\def\basic{1} \def\tbasic{Basic theory}

\def\examplesrepresentations{\basic.1}
\def\groupalgebra{\basic.2}
\def\subrepresentationalgebra{\basic.3}
\def\characterizationirreducible{\basic.4}
\def\decompositionnotunique{\basic.5}
\def\SchurI{\basic.6}
\def\anotherSchur{\basic.7}
\def\decomposition{\basic.8}
\def\standardirreducible{\basic.9}
\def\fixedrank{\basic.10}
\def\Brion{\basic.11}
\def\irreduciblesurjective{\basic.12}
\def\firstdecomposition{\basic.13}

\def\idempotents{2} \def\tidempotents{Decomposition of $K(G)$ via idempotent elements}

\def\canonicaldecomposition{\idempotents.1}
\def\idempotentsdecomposition{\idempotents.2}
\def\characterizationprimitive{\idempotents.3}
\def\primitiverepresentation{\idempotents.4}
\def\propertiescenter{\idempotents.5}
\def\exampleSthree{\idempotents.6}
\def\multiplepoints{\idempotents.7}
\def\examplecyclic{\idempotents.8}
\def\idealcenter{\idempotents.9}
\def\evaluationisomorphism{\idempotents.10}
\def\nonilpotents{\idempotents.11}
\def\ArtinWedderburn{\idempotents.12}

\def\field{3} \def\tfield{Splitting fields}

\def\findrepresentation{\field.1}
\def\SchurSthree{\field.2}
\def\charactersidempotents{\field.3}
\def\quaternions{\field.4}
\def\dihedral{\field.5}
\def\mainconstruction{\field.6}
\def\propertiessplitting{\field.7}
\def\characterspoints{\field.8}
\def\duality{\field.9}
\def\basisevaluation{\field.10}
\def\remarkequations{\field.11}

\def\symmetries{4} \def\tsymmetries{Application to symmetries of functions}

\def\Jsymmetry{\symmetries.1}
\def\idempotentssymmetric{\symmetries.2}
\def\coordinatessymmetric{\symmetries.3}
\def\quadrilinearrepresentations{\symmetries.4}
\def\decompositionofforms{\symmetries.5}
\def\examplessymmetries{\symmetries.6}
\def\orthogonalsymmetric{\symmetries.7}
\def\hyperdeterminant{\symmetries.8}
\def\badidea{\symmetries.9}

\def\Burnside{[B]}
\def\Fulton{[F]} 
\def\FultonHarris{[FH]}
\def\Harris{[H]}

\def\Lam{[L]}
\def\LuxPahlings{[LP]}
\def\MRS{[MRS]}

\def\Quintanilla{[Q]}
\def\Tocino{[T]}

\def\A{\Bbb A}
\def\qed{\hfill\vbox{\hrule\hbox{\vrule\kern3pt
\vbox{\kern6pt}\kern3pt\vrule}\hrule}\bigskip}

\def\mapright#1{\smash{
   \mathop{\longrightarrow}\limits^{#1}}}

\def\mapdown#1{\Big\downarrow
   \rlap{$\vcenter{\hbox{$\scriptstyle#1$}}$}}

\centerline{\bf REPRESENTATION THEORY OF FINITE GROUPS}
\centerline{\bf THROUGH (BASIC) ALGEBRAIC GEOMETRY}
\centerline{Enrique Arrondo}
\bigskip
\bigskip


\noindent{\bf Abstract.} We introduce a new approach to representation theory of finite groups that uses some basic algebraic geometry and allows to do all the theory without using characters. With this approach, to any finite group $G$ we associate a finite number of points from which we can decide which ground fields work fine for the representations of $G$. We apply this point of view to the symmetric group $S_d$, finding easy equations for the different symmetries of functions in $d$ variables. As a byproduct, we give an easy proof of a recent result by Tocino that states that the hyperdeterminant of a $d$-dimensional matrix is zero for all but two types of symmetry.
\bigskip


\bigskip
\bigskip

\centerline{\bf Introduction}
\bigskip

The aim of this note is to present a survey on representation theory of finite groups from a new point of view, based on some basic algebraic geometry. This will allow to recover all the main results of representation theory avoiding the character theory, using instead the --from my point of view-- more natural tool of idempotent elements. Moreover, this approach also allows to discuss which is the natural ground field to work with for each group $G$. This will depend --but not only-- on the coordinates of a finite set of points naturally defined from the group $G$.

Since we are putting together two different branches of mathematics, we would like this paper to be as much self-contained as possible for anyone who is not familiar with at least one of the branches (as a consequence, experts in one of the fields will find unnecessary and very annoying many of the details). We are using only two results from algebraic geometry, which are not very difficult. The first one, in Lemma \fixedrank, is the dimension of a determinantal variety; however we were suggested an alternative proof based only on basic linear algebra. The second one, in Remark \idealcenter, would typically require the Nullstellensatz. However, in our particular case of a finite number of points, we give an alternative easy proof (see Proposition \evaluationisomorphism). About representation theory, we will recall all the basic notions and prove all the results, although many of the proofs follow our new point of view (the standard proofs we give can be found in any basic text on representation theory, as \FultonHarris). Part of this has been done by Patricia Quintanilla in her Bachelor Thesis \Quintanilla, presented under my supervision.

In a first section, we will recall the basic notions of the representation theory of a finite group $G$, such as Schur Lemma, decomposition into irreducible representations or the study of the group algebra $K(G)$, which is the main object to study. To get the nice expected results in representation theory, for simplicity we will assume many times in this first section, that the ground field $K$ is algebraically closed. Indeed, as the experts in the field know, Maschke's theorem states that $K(G)$ is semisimple when the characteristic of $K$ does not divide $|G|$, and thus Artin-Wedderburn theorem implies thus that, when $K$ is good enough, $K(G)$ is the direct product of finite number of matrix algebras. In the language of representation theory, this means that the natural action of $G$ on $K(G)$ (the so-called regular representation) can be decomposed as $K(G)\cong\bigoplus_iEnd(V_i)$, where the $V_i$'s are all the irreducible representations of $G$. This decomposition is the goal of the second and third sections, and we will do it via idempotent elements. 

In the second section, we will see that such idempotents can be computed at once from a finite set of point naturally associated with the group $G$. It is essentially here that we will use some algebraic geometry, which we will present to non-experts in the field.

In a third section we will show that we will need more than the set of points if we want a better behavior, in particular the Artin-Wedderburn theorem. As a byproduct, we will see that our approach with the set of points is naturally equivalent to the approach of the theory of characters (however, our approach has the advantage that the set of points associated with the group indicates a priori which fields $K$ are suitable for representations of the group). We will finish the section with a summary of our main construction.

Finally, in the fourth section we apply the previous techniques to study the decomposition of functions in $d$ variables into their different types of symmetric functions. This corresponds to the study of representations of the symmetric group $S_d$. Finally, we show that our approach allows to find equations for the symmetries much simpler than the ones given in \MRS. Using these simpler equations we gave an easier proof of the result of Tocino (\Tocino) that the hyperdeterminant of a $d$-dimensional matrix is zero for all types of symmetries but two.

\noindent{\bf Acknowledgments.} This research was developed in the framework of the projects MTM2015-65968-P and PID2021-124440NB-I00 funded by the Spanish Government. I thank Lucio Centrone, who encouraged me to write down this paper. I have to thank Hiraku Atobe, who communicated me a second proof of the good version of Lemma \fixedrank, which avoids algebraic geometry. I also thank Maarten Solleveld, Robert Guralnick, Michel Brion and J\"urgen M\"uller for pointing me that I needed some extra hypothesis in the proof of the previous version of Lemma \fixedrank\ and Theorem \irreduciblesurjective. In fact, a fruitful interplay with J\"urgen M\"uller gave me a much better understanding of the topic and helped me to, hopefully, present the contents of this note in a much clearer way. Also Giorgio Ottaviani gave me several useful hints to improve the presentation.

\bigskip
\bigskip
\noindent{\bf \S\basic. \tbasic.} 
\medskip

\noindent{\bf Definition.} A {\it representation of a group} $G$ on a vector space $V$ over a field $K$ is a group homomorphism $\rho:G\to Aut(V)$. We will usually write $gv:=(\rho(g))(v)$. Given two representations $\rho:G\to Aut(V)$ and $\rho':G\to Aut(V')$ a linear map $f:V\to V'$ is called a {\it morphism of representations} if $f(gv)=gf(v)$ for all $g\in G$ and all $v\in V$. Equivalently, the linear map $f$ satisfies $f\circ\rho(g)=\rho'(g)\circ f$ for all $g\in G$

In other words, a representation is a way of considering the group as a subgroup of matrices with entries in $K$. One could wonder about the linear span of those matrices. In the case of finite groups (this will be the case in this note), an intrinsic way of doing so is to consider the group algebra $K(G)$, the vector space with basis $\{e_g\ |\ g\in G\}$; the algebra structure is given by $e_{g'}e_g=e_{g'g}$. In this way, we can extend $\rho$ to a homomorphism of $K$-algebras $\tilde\rho:K(G)\to End(V)$, in which the image of $\alpha=\sum_ga_ge_g\in K(G)$ is the endomorphism $\tilde\rho(\alpha)$ mapping each $v\in V$ to $\alpha\cdot v:=\sum_ga_g(gv)$.

\noindent{\bf Example \examplesrepresentations.} We give a first list of examples that are very simple, but they will be relevant for our purposes.

1) The {\it trivial representation} of any group $G$ is the action of $G$ on $K$ defined by $g\lambda=\lambda$ for any $g\in G$ and any $\lambda\in K$. If $V$ is an $n$-dimensional vector space and we define the trivial action $gv=v$ for any $g\in G$ and any $v\in V$, this representation is isomorphic to the direct sum $n$ times of the trivial representation. Indeed, fixing a basis $v_1,\dots,v_n$ of $V$, the map $K^n\to V$ mapping any $(\lambda_1,\dots,\lambda_n)\in K^n$ to $\lambda_1v_1+\dots+\lambda_nv_n$ is an isomorphism of representations.

2) The {\it regular representation} of any group $G$ is the natural action of $G$ on $K(G)$ given by $g\alpha:=e_g\alpha$. As we will see next, if we want $\tilde\rho:K(G)\to End(V)$ to be a morphism of representations, then we need to endow $End(V)$ with a particular representation structure. Observe that, in a more general framework, if $V,W$ are two representations, the standard way of giving $Hom(V,W)$ a representation structure is by letting $g\in G$ act on $f:V\to W$ as the linear map $gf$ mapping each $v$ to $gf(g^{-1}v)$. When $V=W$, this is different from the structure given in in Lemma \groupalgebra.

3) If $G=S_d$ is the group of permutations of $d$ elements, the {\it alternating representation} is the action of $S_d$ on $K$ given by $\sigma\lambda=sgn(\sigma)\lambda$ for any permutation $\sigma\in S_d$ and any $\lambda\in K$.

We start with a series of basic results.

\proclaim Lemma \groupalgebra. For each representation $\rho:G\to Aut(V)$ on a $K$-vector space $V$, let $\tilde\rho:K(G)\to End(V)$ its corresponding homomorphism of $K$-algebras. 
\item{(i)} If we let $G$ act on $End(V)$ by $gf=\rho(g)\circ f$, then $\tilde\rho$ is a morphism of representations.
\item{(ii)} Fixing a basis $v_1,\dots,v_n$ of $V$ (hence $V$ has dimension $n$), the map $\varphi:End(V)\to V^n$ mapping each $f$ to $\big(f(v_1),\dots,f(v_n)\big)$ is an isomorphism of representations.
\item{(iii)} For each $\alpha\in K(G)$, let $\alpha\cdot:K(G)\to K(G)$ be the left multiplication by $\alpha$ and let $\tilde\rho(\alpha)\circ:End(V)\to End(V)$ be the map defined by the left composition with $\tilde\rho(\alpha)$. Then the diagram
$$\matrix{
K(G)&\mapright{\alpha\cdot}&K(G)\cr
\mapdown{\tilde\rho}&&\mapdown{\tilde\rho}\cr
End(V)&\mapright{\tilde\rho(\alpha)\circ}&End(V)
}$$
is commutative.
\item{(iv)} In the conditions of (ii) and (iii), if $(\alpha,\dots,\alpha)$ denotes the diagonal map $V^n\to V^n$ in which each diagonal map $V\to V$ is the left multiplication by $\alpha$, the diagram
$$\matrix{
End(V)&\mapright{\tilde\rho(\alpha)\circ}&End(V)\cr
\mapdown{\varphi}&&\mapdown{\varphi}\cr
V^n&\mapright{(\alpha,\dots,\alpha)}&V^n
}$$
is commutative.

\noindent {\it Proof.} Part (i) is obvious, since 
$$\tilde\rho(g\alpha)=\tilde\rho(e_g\alpha)=\tilde\rho(e_g)\circ\tilde\rho(\alpha)=\rho(g)\circ\tilde\rho(\alpha)=g\tilde\rho(\alpha).$$

For part (ii), it is clear that $\varphi$ is an isomorphism of vector spaces, so it suffices to prove that it is a morphism of representations. This is so because, for each $g\in G$ and $f\in End(V)$, one has
$$\varphi(gf)=\big((gf)(v_1),\dots,(gf)(v_n)\big)=\big(gf(v_1),\dots,gf(v_1)\big)=g\big(f(v_1),\dots,f(v_n)\big)=g\varphi(f).$$

For part (iii), we easily see that, for each $\beta\in K(G)$ and any $v\in V$, one has $\tilde\rho(\alpha\beta)(v)=(\alpha\beta)v=\alpha(\beta v)=\big(\tilde\rho(\alpha)\big)(\beta v)=\big(\tilde\rho(\alpha)\big)\big(\tilde\rho(\beta)(v)\big)$.

Finally, part (iv) follows because, for each $f\in End(V)$, one has
$$\varphi\big(\tilde\rho(\alpha)\circ f\big)=\big(\alpha f(v_1),\dots,\alpha f(v_n)\big)=(\alpha,\dots,\alpha)\big(f(v_1),\dots,f(v_n)\big)=(\alpha,\dots,\alpha)(\varphi(f)).$$
\qed

The following trivial remark will be key to the whole theory:

\proclaim Lemma \subrepresentationalgebra. Any subrepresentation $W\subset K(G)$ of the regular representation has a natural structure of subalgebra.

\noindent{\it Proof.} It is enough to prove that the product in $K(G)$ of two elements $w_1,w_2\in W$ is still in $W$. This is obvious, since $w_1$ takes the form $w_1=\lambda_1e_{g_1}+\dots+\lambda_re_{g_r}$, with $\lambda_1,\dots,\lambda_r\in K$ and $g_1,\dots,g_r\in G$, and each $e_{g_i}w_2=g_iw_2$ is in $W$.
\qed

\noindent{\bf Definition.} An {\it irreducible representation} is a representation not possessing non-trivial subrepresentations.

When the ground field has not a bad characteristic, the notion of irreducibility is the expected one:

\proclaim Theorem \characterizationirreducible. Assume that $|G|$ is not divisible by the characteristic of $K$. Then, for any subrepresentation $V'$ of a representation $V$ there exists another subrepresentation $V''$ such that $V=V'\oplus V''$. As a consequence, a representation is irreducible if and only if it does not split in a nontrivial way as a direct sum of two representations.

\noindent{\it Proof.} Take any complementary subspace $W$ of $V'$ in $V$. If $p_1:V=V'\oplus W\to V'$ is the first projection, we define $f:V\to V'$ as
$$f(v)={\sum_{g\in G}gp_1(g^{-1}v)\over|G|}$$
which is a morphism of representations restricting to the identity on $V'$ (in particular, $f\circ f=f$). Hence its kernel $V''$ is a subrepresentation, and, since any $v\in V$ can be written as $v=f(v)+(v-f(v))$, it follows $V=V'\oplus V''$.
\qed

\noindent{\bf Remark \decompositionnotunique.} The complementary subspace $W$ of the above result is not unique. For example, the isomorphism of Lemma \groupalgebra(ii) induces a decomposition $End(V)=V'_1\oplus\dots\oplus V'_n$, where $V'_i$ is the subspace of endomorphisms vanishing at the hyperplane $H_i\subset V$ generated by $v_1,\dots,v_{i-1},v_{i+1},\dots,v_n$. Hence, the subrepresentation $V'_n$ has a complementary subrepresentation $V'_1\oplus\dots\oplus V'_{n-1}$ consisting of the endomorphisms vanishing at $v_n$. Obviously, changing this last vector of the basis by another one that is not multiple of $v_n$, produces a different complementary subrepresentation. Although the decomposition of $End(V)$ is not unique, observe that any of the possible decompositions we found comes from an isomorphism of representations $End(V)\cong V^{\oplus \dim(V)}$. This is the kind of uniqueness that we will prove in Theorem \decomposition. 

\proclaim Proposition \SchurI\ {\rm (Schur's Lemma)}. A nonzero morphism between irreducible representations is necessarily an isomorphism.

\noindent{\it Proof.} Let $f:V\to W$ be a nonzero morphism of irreducible representations. Since $\ker(f)$ is a proper subrepresentation of $V$, necessarily is the zero representation, so that $f$ is injective. Also Im$(f)$ is a nonzero subrepresentation of $W$, so that it is $W$ and hence $f$ is surjective.
\qed

\noindent{\bf Remark \anotherSchur.} Schur's Lemma has a second part, proved easily when $K$ is the complex field (or any algebraically closed field). We will obtain it in general in Theorem \propertiessplitting(ii). 
 
\proclaim Theorem \decomposition. If $|G|$ is not divisible by the characteristic of $K$, any representation of $G$ decomposes as a direct sum of irreducible representations. Such a decomposition is, up to isomorphism, unique. 

\noindent{\it Proof.} We will prove the result by induction on the dimension of the representation, the result being trivial if the dimension is one. Let $V$ be a representation of $G$ and let us see that it decomposes as a direct sum of irreducible representations. If $V$ is irreducible, there is nothing to prove. If it is not irreducible, Theorem \characterizationirreducible\ implies that we can write $V=V'\oplus V''$ in a nontrivial way. Hence, by induction hypothesis, $V'$ and $V''$ decompose as a direct sum of irreducible representations. The uniqueness is an easy consequence of Schur's Lemma.
\qed

\noindent{\bf Example \standardirreducible.} Consider the natural representation of the symmetric group $S_d$ on the vector space $K^d$ with canonical basis $e_1,\dots,e_d$ by defining, for each permutation $\sigma$, $\sigma e_i=e_{\sigma(i)}$. This representation is not irreducible, since the span of $e_1+\dots+e_d$ and the hyperplane $V\subset K^d$ of equation $x_1+\dots+x_d=0$  are invariant under the action. If the characteristic of $K$ is a divisor of  $d$, the span of $e_1+\dots+e_d$ is again an invariant subspace of $V$, so that $V$ is not irreducible for $d\ge3$. If instead $d$ is not divisible by the characteristic of $K$, $K^d$ decomposes as the direct sum of the span of $e_1+\dots+e_d$ (which is isomorphic to the trivial representation) and $V$ (this is usually called the {\it standard representation} of $S_d$). Let us see that, in this case, $V$ is irreducible.

\noindent To prove this, we use the following trivial observation: if a representation $\rho:G\to Aut(V)$ contains an invariant subrepresentation $W\subset V$, obviously the image of $\tilde\rho$ is contained in the proper subspace of $End(V)$ consisting of the endomorphisms preserving $W$, so that $\tilde\rho$ is not surjective. Hence, the surjectivity of $\tilde\rho$ implies the irreducibility of the representation. In our example, we first consider the basis of $V$ given by $e_1-e_d,\dots,e_{d-1}-e_d$. Thus a basis of $End(V)$ is given by the endomorphisms $\varphi_{ij}$ mapping $e_i-e_d$ to $e_j-e_d$ and any other $e_k-e_d$ of the basis to zero. It is enough to prove that each $\varphi_{ij}$ is in the image of $\tilde\rho$. For this, we first consider $\Sigma_{l\ne i}e_{(l\ i)}-(d-2)e_{(1)}$, which maps each $e_k-e_d$ to zero, except $e_i-e_d$, which is mapped to $\Sigma_{l\ne i}e_l-(d-1)e_i$. The map $e_{(i\ d)}-e_{(1)}$ sends the latter to $d(e_i-e_d)$. As a consequence, 
the product ${1\over d}e_{(i\ j)}(e_{(i\ d)}-e_{(1)})(\Sigma_{l\ne i}e_{(l\ i)}-(d-2)e_{(1)})$ maps to $\varphi_{ij}$, proving the surjectivity of $\tilde\rho$.

\bigskip

It is an important result, due to Burnside (\Burnside) in characteristic zero (see \Lam\ for arbitrary characteristic), that the surjectivity of $\tilde\rho$ in the previous example is in fact a characterization of the irreducibility of a presentation when the field $K$ is sufficiently good (this is the notion of splitting field that we will explain in section \field). Let us see that it is easy to prove this characterization first when the field is algebraically closed, and one of the main goals of the paper is to see for which other fields a result like that still remains true.

The result we will need is the following, which can be obtained either from algebraic geometry, although we will include another proof using elementary linear algebra, indicated to me by Hiraku Atobe.

\proclaim Lemma \fixedrank. Let $V',V$ be vector spaces of respective dimensions $0<n'<n$ over an algebraically closed field $K$. Assume $Hom(V',V)$ has a linear subspace $W$ of dimension $n$ such that any nonzero $f\in W$ has rank $n'$. Then $n'=1$.

\noindent {\it Proof.} Using algebraic geometry, it is very well-known (see for Example \Harris\ Proposition 12.2) that the set of homomorphisms of rank at most $n'-1$ has codimension $n-n'+1$ in $Hom(V',V)$. Hence if it were $n'>1$ the subspace $W$ would necessarily contain a nonzero homomorphism of rank at most $n'-1$, contradicting our hypothesis.

Without algebraic geometry, observe first that, observe that any nonzero $f\in W$ is injective, because it has maximal rank. This implies that, for any nonzero $v'\in V$, the evaluation map $ev_{v'}:W\to V$ defined by $ev_{v'}(f)=f(v')$ is injective, hence an isomorphism (because $W$ and $V$ has the same dimension). Now, given any two nonzero $v',w'\in V'$, the composition 
$$ev_{w'}\circ ev_{v'}^{-1}:V\to W\to V$$
is an endomorhism which necessarily (since $K$ is algebraically closed) has a non zero eigenvector $v$ with eigenvalue $\lambda$. This means that $f:=ev_{v'}^{-1}(v)$ satisfies $f(w')=\lambda f(v')$. The injectivity of $f$ implies $w'=\lambda v'$, concluding that $V'$ has necessarily dimension one.
\qed

\noindent{\bf Remark \Brion.} The hypothesis that $K$ is algebraically closed is crucial. A simple example, suggested to me by Michel Brion, is the following. Take a chain of finite extensions $K\subset L'\subset L$. Let $n',n$ be respectively the dimensions of $L',L$ as vector spaces over $K$. We consider the subspace $W\subset Hom(L',L)$ consisting of the maps $\varphi_x:L'\to L$ (when $x\in L$) defined by $\varphi_x(y)=xy$. Since any $\varphi_x$ is injective for any $x\ne0$, they have rank $n'$. Clearly $W$ has dimension $n$, because the map $L\to W$ mapping $x$ to $\varphi_x$ is an isomorphism. Taking any $n'>1$ we get a counterexample to the above Lemma.

We can now prove the promised characterization:

\proclaim Theorem \irreduciblesurjective. If $K$ is algebraically closed, a representation $\rho:G\to Aut(V)$ is irreducible if and only if the map $\tilde\rho$ is surjective. 

\noindent {\it Proof.} We only need to proof that the surjectivity is a necessary condition. The main tool is that, for each nonzero vectors $v,v'\in V$, there is $\alpha\in K(G)$ such that $\alpha v=v'$. To prove this we need to prove that the map $\tilde\rho_v:K(G)\to V$ sending each $\alpha$ to $\alpha v$ is surjective. This is so because that map is the composition of $\tilde\rho$ with the evaluation map $End(V)\to V$ at $v$, and both maps are morphisms of representations. Since $V$ is irreducible and the map is not zero (because the image of $1$ is $v$), the map is surjective.

Now fix a nonzero $\alpha\in K(G)$ such that $\tilde\rho(\alpha)$ has minimal rank $r$, and fix a basis $v_1,\dots,v_n$ of $V$. By the surjectivity of $\tilde\rho_{\alpha v}$, we have $\alpha_1,\dots,\alpha_n\in K(G)$ such that $\alpha_i(\alpha v)=v_i$ for $i=1,\dots,n$. This implies that the endomorphisms $\tilde\rho(\alpha_1\alpha),\dots,\tilde\rho(\alpha_n\alpha)$ are linearly independent. Since any nonzero linear combination of them clearly has rank at most $r$ (its kernel contains the kernel of $\tilde\rho(\alpha)$), the minimality of $r$ and Lemma \fixedrank\ imply $r=1$ (by regarding each $\tilde\rho(\alpha_i\alpha)$ as a linear map from $V/\ker(\tilde\rho(\alpha))\to V$). 

At this point we construct a suitable basis for $V$. We start by taking $v_1$ a vector outside $H_1:=\ker(\tilde\rho(\alpha))$, which we just proved to be a hyperplane of $V$, and write $f_1:=\tilde\rho(\alpha)$. We then take a nonzero vector $v_2\in H_1$, and hence there exists $\alpha_2\in K(G)$ such that $\alpha_2v_2=v_1$. This implies that $f_2:=\tilde\rho(\alpha\alpha_2)$ is an endomorphism of rank one, such that $v_2$ is not in its kernel $H_2$. In particular, $H_1\not\subset H_2$, so that $H_1\cap H_2$ has codimension two. We can iterate the process to construct $v_1,v_2,\dots,v_n$ such that, for each $i$, there is an endomorphism $f_i:V\to V$ in the image of $\tilde\rho$ whose kernel is a hyperplane $H_i\supset v_{i+1},\dots,v_n$ and $f_i(v_i)=\alpha v_1\ne0$. This can be done picking, for each $i$, a nonzero vector $v_i\in H_1\cap\dots\cap H_{i-1}$, an element $\alpha_i\in K(G)$ such that $\alpha_iv_i=v_1$ and $f_i=\tilde\rho(\alpha\alpha_i)$. 

It is now straightforward to check that $v_1,\dots,v_n$ is a basis. Moreover, if for each $i,j$ we consider the composition $f_{ij}$ of the above $f_i$ with an endomorphism in the image of $\tilde\rho$ mapping $f_i(v_i)$ to $v_j$, the matrix of $f_{ij}$ with respecto to the basis $v_1,\dots,v_n$ takes the form
$$\pmatrix{
a_{11}&\dots&a_{i-1,1}&0&0&\dots&0\cr
\vdots&&\vdots&\vdots&\vdots&&\vdots\cr
a_{1,j-1}&\dots&a_{i-1,j-1}&0&0&\dots&0\cr
a_{1j}&\dots&a_{i-1,j}&1&0&\dots&0\cr
a_{1,j+1}&\dots&a_{i-1,j+1}&0&0&\dots&0\cr
\vdots&&\vdots&\vdots&\vdots&&\vdots\cr
a_{1n}&\dots&a_{i-1,n}&0&0&\dots&0\cr}.$$
It is easy to check that these matrices form a basis of the space of $n\times n$ matrices, hence the endomorphisms $f_{ij}$ also form a basis of $End(V)$. Since all them are in the image of $\tilde\rho$, this proves the wanted surjectivity.
\qed

We finish with a key decomposition in the general case, in which we will not assume a priory the surjectivity of $\tilde\rho$, but that we will inspire the main technique in the next section.

\proclaim Proposition \firstdecomposition. If the characteristic of $K$ does not divide $|G|$, for any irreducible representation $\rho:G\to Aut(V)$, there is a decomposition $K(G)=W\oplus\ker(\tilde\rho)$ as unitary algebras, where $W$ is isomorphic to the image of $\tilde\rho:K(G)\to End(V)$. 

\noindent {\it Proof.} The condition about the characteristic implies, by Theorem \characterizationirreducible, that $K(G)$ decomposes, as a representation, as $W\oplus\ker(\tilde\rho)$, and $W$ is isomorphic to the image of $\tilde\rho:K(G)\to End(V)$. By Lemma \subrepresentationalgebra, this also holds as algebras. Since the image of $\tilde\rho$ contains $id_V$, also $W$ is unitary. Hence, the difference of $1\in K(G)$ and the unit in $W$ is also a unit in $\ker(\tilde\rho)$.
\qed

\bigskip
\noindent{\bf \S\idempotents. \tidempotents.} 
\medskip

To understand the main idea of this section, we start with the following clarifying remark:

\noindent{\bf Remark \canonicaldecomposition.} In view of Remark \decompositionnotunique, one could guess that the decomposition in Proposition \firstdecomposition\ is not unique. However, this is not the case now. The reason is that, considering the corresponding unit elements $\beta\in W$ and $\beta'\in\ker(\tilde\rho)$, they are in the center of $K(G)$,  their sum is the unit element of $K(G)$ and clearly $W=\beta K(G)$ and $\ker(\tilde\rho)=\beta'K(G)$. In particular, the subspace $W$ is unique, since depends on $\beta=1-\beta'$, where $\beta'$ is the unit element of the kernel of $\tilde\rho$, which in turn depends on the irreducible representation $V$. 

Hence the decomposition of $K(G)$ comes from the existence of {\it central idempotents}, i.e. idempotent elements in the center of the algebra. This is what we are going to study in this section, and we will get from the idempotent elements a first decomposition which will be very close to the wanted Artin-Wedderburn decomposition (in the next section we will see what is missing to get the right decomposition).

To start with, the precise keystone result, in the more general framework of rings, is the following:

\proclaim Lemma \idempotentsdecomposition. A unitary ring $A$ decomposes as a direct sum of two non-trivial unitary rings if and only if $A$ possesses a central idempotent $\beta\ne0,1$. Moreover, in such a case, $1-\beta$ is also a central idempotent and one has the decomposition $A=\beta A\oplus(1-\beta)A$, the respective unit elements being $\beta$ and $1-\beta$.

\noindent{\it Proof.} If there exists a non-trivial decomposition $A=B\oplus C$ as unitary rings with respective units $\beta,\gamma\ne0,1$, it is clear that $\beta$ and $\gamma$ are idempotent elements in the center of $A$, and $1=\beta+\gamma$.

Reciprocally, given an idempotent element $\beta\ne0,1$ in the center of $A$, it is straightforward to check that $1-\beta$ is also an idempotent element in the center of $A$ and one has the decomposition $A=\beta A\oplus(1-\beta)A$, the respective unit elements being $\beta$ and $1-\beta$. For example, assume that we have an element $\beta a=(1-\beta)a'$ in the intersection of $\beta A$ and $(1-\beta)A$. Then:
$$\beta a=\beta^2a=\beta(1-\beta)a'=(\beta-\beta^2)a'=0\cdot a'=0.$$
\qed

We characterize next when we cannot decompose the above pieces. 

\proclaim Lemma \characterizationprimitive. Let $A$ be a unitary ring and let $\beta\in A$ be a central idempotent. Then the following are equivalent.
\item{(i)} $\beta$ does not decompose as a sum of two non-trivial central idempotents.
\item{(ii)} $\beta A$ does not decompose as a direct sum of non-trivial unitary rings.
\item{(iii)} The center of $\beta A$ does not decompose as a direct sum of non-trivial unitary rings.

\noindent{\it Proof.} From Lemma \idempotentsdecomposition\ we get that (i) is equivalent to (ii), because $\beta$ is the unit element of $\beta A$. On the other hand, a nontrivial decomposition of $\beta A$ yields a non trivial decomposition of its center, since the center of each piece necessarily contains its own unit element; hence (iii) implies (ii). Finally, by Lemma \idempotentsdecomposition, a decomposition of the center of $\beta A$ implies the decomposition of its unit element, which is $\beta$, so that (i) implies (iii).
\qed

\noindent{\bf Definition.} A {\it centrally primitive idempotent} of a ring $A$  is a central idempotent satisfying any of the equivalent conditions of Lemma \characterizationprimitive.

\proclaim Theorem \primitiverepresentation. If $|G|$ is not divisible by the characteristic of $K$, for any centrally primitive idempotent $\beta\in K(G)$, there is a unique irreducible representation $\rho:G\to Aut(V)$ such that $\beta K(G)$ is isomorphic to the image of $\tilde\rho:K(G)\to End(V)$. 

\noindent{\it Proof.} Let $\beta$ be a centrally primitive idempotent, and let $V\subset\beta K(G)$ be an irreducible subrepresentation. Since $(1-\beta)\beta K(G)=0$, then also $(1-\beta)V=0$, so that $\beta V=V$. In particular (see Remark \decompositionnotunique), $\beta End(V)=End(V)$ and $\beta W'=W'$ where $W'$ is the image of $\tilde\rho:K(G)\to End(V)$.  By Proposition \firstdecomposition, we can decompose $K(G)=W\oplus\ker\tilde\rho$, with $W$ isomorphic to $W'$, so that $\beta K(G)=\beta W\oplus\beta\ker\tilde\rho=W\oplus\beta\ker\tilde\rho$. Since $\beta$ is primitive, necessarily $\beta\ker\tilde\rho=0$ and hence $\beta K(G)=W$. Since $W$ is a subrepresentation of $End(V)$, this implies (see Remark \decompositionnotunique) that the decomposition of $\beta K(G)$ into irreducible representations is the direct sum of copies of $V$. The unicity of the decomposition implies that $V$ is unique.
\qed

In order to completely decompose a ring $A$, one should find its centrally primitive idempotent elements. And, as observed in Lemma \characterizationprimitive, it is enough to decompose its center $ZA$. In the case of $K(G)$, its center $ZK(G)$ is a commutative algebra which is {\it finite}, i.e. finite dimensional as a $K$-vector space. Of course, we could be sure that the decomposition is complete if the pieces are one-dimensional. This is exactly what we are going to prove. We first study the center. Before doing this, we will fix some notation that we will use throughout the paper:

\noindent {\bf Notation.} We will write $P$ for the set of conjugacy classes of a finite group $G$. When $G=S_d$, we identify its set of conjugacy classes with the set $P_d$ of partitions of $d$. This is so because two elements of $S_d$ are in the same conjugacy class if and only if they both decompose into a product of disjoint cycles of the same lengths $\lambda_1\ge\dots\ge\lambda_r$ (we include here cycles of length one for fixed elements of the permutation), in which obviously $d=\lambda_1+\dots+\lambda_r$ (for example, $\lambda=(1,\dots,1)$ is the conjugacy class of the identity). We will write as $\lambda=(\lambda_1,\dots,\lambda_r)$ for such a partition. In analogy with this, we will typically use $\lambda$ to denote a conjugacy class of an arbitrary group $G$. We will also write $P'$ for the set of conjugacy classes different from the class of the unit element.

\proclaim Proposition \propertiescenter. Let $G$ be a finite group and let $ZK(G)$ be the center of the group algebra $K(G)$. Then
\item{(i)} For any $\alpha\in ZK(G)$ and any representation $V$, the left multiplication by $\alpha$ is a morphism of representations $\alpha\cdot:V\to V$.
\item{(ii)}  The center of $K(G)$ has, as a vector space, a basis consisting of the set of elements
$$\alpha_\lambda:=\sum_{g\in\lambda}e_g$$
where $\lambda$ varies in the set $P$ of conjugacy classes in $G$. 

\noindent{\it Proof.} Part (i) is immediate because, by definition of morphism of representations, the multiplication by $\alpha$ must commute with any $\rho(g)$, i.e. $\alpha$ must commute with any $e_g$, and this is the case precisely because $\alpha$ is in the center of $K(G)$. 

For part (ii), let $\alpha=\sum_{g\in G}\mu_ge_g$ be an element of the center of $K(G)$. This is equivalent to say that $\alpha e_h=e_h\alpha$ for all $h\in G$, i.e. $\sum_{g\in G}\mu_ge_{gh}=\sum_{g\in G}\mu_ge_{hg}$. Since the coefficient of $e_{gh}$ in the right-hand sum is $\mu_{h^{-1}gh}$, it follows that $\alpha$ is in $ZK(G)$ if and only if, for all $g,h\in G$, $\mu_g=\mu_{h^{-1}gh}$. In other words, the coefficient of $e_g$ in $\alpha$ coincides with the coefficient of any conjugate $e_{h^{-1}gh}$. This means that $\alpha$ is a linear combination of the elements $\alpha_\lambda$, as wanted.
\qed

\noindent{\bf Example \exampleSthree.} Let us illustrate what we want to do in the simple example $G=S_3$, the group of permutations of three elements. According to the above notation with partitions, the conjugacy classes of $S_3$ are:

--The class of the identity, which is, as a product of disjoint cycles, $id=(1)(2)(3)$, so that it corresponds to the partition $\lambda=(1,1,1)$.

--The class of the $2$-cycles (hence a product of a transposition and a $1$-cycle), which corresponds to the partition $\lambda=(2,1)$.

--The class of the $3$-cycles, i.e. correspondinf to the partition $\lambda=(3)$.
 
 \noindent As stated in Proposition \propertiescenter(ii), the finite algebra $ZK(S_3)$ has a basis, as a vector space, given by
$$1:=\alpha_{(1,1,1)}=e_{(1)}$$
$$\alpha_{(2,1)}:=e_{(1\ 2)}+e_{(1\ 3)}+e_{(2\ 3)}$$
$$\alpha_{(3)}:=e_{(1\ 2\ 3)}+e_{(1\ 3\ 2)}$$
and the algebra structure is given by the relations
$$\alpha_{(2,1)}^2=3+3\alpha_{(3)}$$
$$\alpha_{(2,1)}\alpha_{(3)}=2\alpha_{(2,1)}$$
$$\alpha_{(3)}^2=2+\alpha_{(3)}.$$
Hence $ZK(S_3)$ is isomorphic to $K[x,y]/I$, where $I$ is the ideal generated by the polynomials $x^2-3y-3,xy-2x,y^2-y-2)$. Solving this system of equations, one get the points $(3,2),(-3,2),(0,-1)$. We can consider the composed map in which the second map is the evaluation of classes of polynomials at these three points
$$\matrix{
ZK(S_3)&\to&K[x,y]/I&\to&K^3\cr
1&\mapsto&\bar 1&\mapsto&(1,1,1)\cr
\alpha_{(2,1)}&\mapsto&\bar x&\mapsto&(3,-3,0)\cr
\alpha_{(3)}&\mapsto&\bar y&\mapsto&(2,2,-1)
}$$
In coordinates, writing $A:=\pmatrix{1&1&1\cr3&-3&0\cr2&2&-1}$,
$$\pmatrix{\mu_1&\mu_2&\mu_3}\pmatrix{1\cr\alpha_{(2,1)}\cr\alpha_{(3)}}\mapsto\pmatrix{\mu_1&\mu_2&\mu_3}A$$
i.e. $A$ is the matrix whose columns are the coordinates of the three points (preceded by a $1$). Since $A$ is a regular matrix, we have an isomorphism, and its inverse is given by $B=A^{-1}$, more precisely
$$\pmatrix{\nu_1&\nu_2&\nu_3}\mapsto\pmatrix{\nu_1&\nu_2&\nu_3}B\pmatrix{1\cr\alpha_{(2,1)}\cr\alpha_{(3)}}=\pmatrix{\nu_1&\nu_2&\nu_3}\pmatrix{1\over6&1\over6&1\over6\cr\cr1\over6&-1\over6&1\over6\cr\cr4\over6&0&-2\over6}\pmatrix{1\cr\alpha_{(2,1)}\cr\alpha_{(3)}}$$
Since the canonical basis of $K^3$ is the set of centrally primitive idempotents and
$$(1,0,0)\mapsto{1+\alpha_{(2,1)}+\alpha_{(3)}\over6}=:\beta_1$$
$$(0,1,0)\mapsto{1-\alpha_{(2,1)}+\alpha_{(3)}\over6}=:\beta_2$$
$$(0,0,1)\mapsto{4-2\alpha_{(3)}\over6}=:\beta_3$$
it follows that $\beta_1,\beta_2,\beta_3$ (whose coordinates are the rows of $B$) are the idempotent elements that allow to decompose $K(S_3)$.

\noindent{\bf Remark \multiplepoints.} In the above example it is crucial that the characteristic of $K$ is not $2$ nor $3$. Observe that, for those particular values of the characteristic, the three points are not different. In algebraic terms, we have that, in characteristic $2$, the ideal $I$ has as primary decomposition
$$I=(x^2+y-1,xy,y^2-y)=(x,y-1)\cap\big(y,(x-1)^2\big)$$
which represents the points $(0,1)$ (counted with multiplicity one) and $(1,0)$ (counted with multiplicity two), while, in characterisitic $3$ the ideal $I=\big(x^2,x(y+1),(y+1)^2\big)$ is primary itself, representing the point $(0,-1)$ with multiplicity three.

\noindent{\bf Example \examplecyclic} Let us study now the case of the cyclic group of orden $n$. To keep the multiplicative notation, we will identify it with $G=\{\omega^i\ |\ i\in{\Bbb Z}\}$ where $\omega\in{\Bbb C}$ is an $n$-th primitive root of $1$. Since $G$ is commutative, $ZK(G)$ coincides with $K(G)$, hence writing $\alpha_i=e_{\omega^i}$,  a basis as a vector is given by $1=\alpha_0,\alpha_1,\dots,\alpha_{n-1}$. Since the relations are $\alpha_i\alpha_j=\alpha_{i+j}$, the commutative algebra $ZK(G)=K(G)$ can be identified with $K[x_1,\dots,x_{n-1}]/I$, where $I$ is the ideal generated by the elements of the form $x_ix_j=x_k$, when $k\equiv i+j\ {\rm mod}\ n$ and the convention $x_0=1$. In this case, the ideal $I$ defines the set of points $a_i=(\omega^i,\omega^{2i},\dots,\omega^{(n-1)i})$, with $i=1,\dots,n$. Observe that, now, if we want the points $a_1,\dots,a_n$ to have their coordinates in $K$, we need $K$ to contain all the $n$-th roots of $1$. This is our goal: given any finite group $G$, to compute its associated set of points to know which fields contain the coordinates of the points. We will see that, to make representation theory on the given group, we will need to restrict to such fields (in general, we will still need some extra conditions, and this will be studied in the next section).

\noindent{\bf Remark \idealcenter.} We generalize now the construction in Examples \exampleSthree\ and \examplecyclic\ to provide the algebra structure of $ZK(G)$ based on the basis of Proposition \propertiescenter(ii). When $\lambda$ is the conjugacy class of $1\in G$, then $\alpha_\lambda=e_1$, which is the unit of the algebras $ZK(G)$ and $K(G)$. If we consider now $P'$ to be the set of conjugacy classes in $G$ different from the class of $1$, there is a natural epimorphism of $K$-algebras
$$\matrix{
K[\{x_\lambda\}_{\lambda\in P'}]&\to&ZK(G)\cr
x_\lambda&\mapsto&\alpha_\lambda
}$$  
Let $I$ be the kernel of that epimorphism. Then $ZK(G)$ can be identified canonically with $K[\{x_\lambda\}_{\lambda\in P'}]/I$. The next step would be, as in Examples \exampleSthree\ and \examplecyclic, to see that the set
$$V(I):=\{a\in\A^{|P'|}_K\ |\ f(a)=0\hbox{ for all }f\in I\}$$
is a finite number of points and that the evaluation at those points decomposes $ZK(G)$ as a direct sum of copies of $K$. To do this, we will need some algebraic geometry again (see, for example, \Fulton \S2.9). The first observation is that, $K[\{x_\lambda\}_{\lambda\in P'}]/I$ being finite-dimensional, it defines necessarily a finite set of points, precisely as many as the dimension of the vector space, but the points could have some multiplicity, as shown in Remark \multiplepoints. The main result is that, by the Nullstellensatz, these points count with multiplicity one precisely when $I$ is a radical ideal. We will prove these facts without the Nullstellesatz in the following:

\proclaim Proposition \evaluationisomorphism. Let $I\subset K[x_1,\dots,x_n]$ be a proper ideal such that $K[x_1,\dots,x_n]/I$ has finite dimension as a vector space over $K$. Then:
\item{(i)} The set $V(I)$ is finite and not empty, and the coordinates of the points of $V(I)$ are in a finite extension of $K$.
\item{(ii)} For any finite set of points $a_1,\dots,a_m\in\A^n_K$, the evaluation map
$$\matrix{
K[x_1,\dots,x_n]&\to&K^m\cr
f&\mapsto&\big(f(a_1),\dots,f(a_m)\big)
}$$
is surjective.
\item{(iii)} If $I$ is radical and the coordinates of the points of $V(I)$ are in $K$, then $I$ is the kernel of the evaluation map at the points $a_1,\dots,a_m$ of $V(I)$.
\item{(iv)} In the conditions of (iii), the number $m$ of points of $V(I)$ coincides with the dimension of $K[x_1,\dots,x_n]/I$, and, as a $K$-algebra, there is an isomorphism $K[x_1,\dots,x_n]/I\cong K^m$ consisting of the evaluation of classes of polynomials at the points of $V(I)$.

\noindent{\it Proof.} For part (i), let $m$ be the dimension of $K[x_1,\dots,x_n]/I$. Then, for each $i=1,\dots,n$, the classes of $1,x_i,x_i^2,\dots,x_i^m$ must be linearly dependent. Hence there exists a nonzero polynomial $f_i(x_i)\in I$, so that the possible coordinates $x_i$ of the points of $V(I)$ are the roots of $f_i$. This implies that $V(I)$ is finite and the coordinates of its points are in any field containing all the roots of $f_1,\dots,f_n$, so that only a finite extension of $K$ is needed. To prove the non-emptyness, since we are considering $V(I)$ inside the affine space over a suitable extension of $K$, we can assume that all the roots of $f_1,\dots,f_n$ are in $K$. Moreover, since any proper ideal is contained in a prime ideal, we can assume $I$ to be prime. Therefore $I$ contains some irreducible factor of each $f_i$. Such factors are linear by assumption, hence there are elements $x_i-a_i\in I$, with $a_i\in K$ for $i=1,\dots,n$. Therefore $I$ contains the maximal ideal $(x_1-a_1,\dots,x_n-a_n)$, thus $I=(x_1-a_1,\dots,x_n-a_n)$, i.e. the ideal of all polynomials vanishing at the point $a=(a_1,\dots,a_n)$ and $V(I)=\{a\}\ne\emptyset$.

For part (ii), it is enough to see that the canonical basis of $K^m$ is in the image of the evaluation map. To prove this, it suffices to find, for each $i=1,\dots,m$, a polynomial $f_i$ vanishing at all $a_1,\dots,a_m$ except $a_i$. For this, it is enough to consider, for each $j\ne i$, the equation $g_j$ of a hypersurface passing through $a_j$ and not through any of the others. We can thus take $f_i=g_1\dots g_{i-1}g_{i+1}\dots g_m$.

For part (iii), the fact that $I$ is radical implies that, in the primary decomposition $I=I_1\cap\dots\cap I_m$ of $I$, all the primary ideals are prime. Reasoning as in part (i), each $I_i$ is the ideal of all polynomials vanishing at some point $a_i\in{\Bbb A}^n_K$. Hence $I$ is the ideal of all polynomials vanishing at $V(I)=\{a_1.\dots,a_m\}$, i.e. $I$ is the kernel of the evaluation map.

Part (iv) is a consequence of (ii) and (iii).
\qed

Since an ideal $I$ of a ring $A$ is radical if and only if $A/I$ has not nonzero nilpotent elements, we need to prove the following:

\proclaim Lemma \nonilpotents. If $|G|$ is not divisible by the characteristic of $K$, the only nilpotent element of $ZK(G)$ is the zero.

\noindent{\it Proof.} Let $\alpha\in ZK(G)$ be a nilpotent element. For any irreducible representation $V$, let us consider the endomorphism $V\to V$ given by the left multiplication by $\alpha$, which is a morphism of representations by Proposition \propertiescenter(i). Since $\alpha$ is nilpotent, the endomorphism is so, hence it cannot have maximal rank, hence by Proposition \SchurI\ is the zero endomorphism. Decomposing $K(G)$ into irreducible representations (for which we need the condition about the characteristic), we conclude that the endomorphism $K(G)\to K(G)$ consisting of the left multiplication by $\alpha$ is zero. Therefore $\alpha=\alpha\cdot e_1=0$, as wanted.
\qed

We can finish the section with a first decomposition result:

\proclaim Theorem \ArtinWedderburn. Let $G$ be a finite group and let $K$ be a field whose characteristic is not a divisor of $|G|$ and containing the coordinates of the points of $V(I)$ defined in Remark \idealcenter. Let $m$ be the number of conjugacy classes of $G$. Then, $G$ possesses exactly $m$ irreducible representations $\rho_1:G\to Aut(V_1)$, $\dots$, $\rho_m:G\to Aut(V_m)$ corresponding to the set $\beta_1,\dots,\beta_m$ of centrally primitive idempotents, and there is a decomposition as algebras $K(G)\cong\bigoplus_{i=1}^mW_i$, where each $W_i=\beta_iK(G)$ is isomorphic to the image of $\tilde\rho_i:K(G)\to End(V_i)$.

\noindent{\it Proof.} By Lemma \nonilpotents, the ideal $I$ of the isomorphism $ZK(G)\cong K[\{x_\lambda\}_{\lambda\in P'}]/I$ is radical. Hence Proposition \evaluationisomorphism\ provides an isomorphism of algebras $ZK(G)\cong K^m$, so that $m$ is the dimension of $ZK(G)$, which is the number of conjugacy classes of $G$ (Proposition \propertiescenter(ii)). Since the unit in $K^m$ is $(1,\dots,1)$, which is the sum of the elements of the canonical basis, the inverse image of the canonical basis of $K^m$ gives a set $\beta_1,\dots,\beta_m$ of centrally primitive idempotents whose sum is $1\in K(G)$. By Lemma \idempotentsdecomposition\ we have a decomposition of unitary algebras
$$K(G)=\beta_1K(G)\oplus\dots\oplus\beta_mK(G)$$
and Theorem \primitiverepresentation\ implies that, for each $i=1,\dots,m$, there is a unique irreducible representation $\rho_i:G\to Aut(V_i)$ such that $\beta_iK(G)$ is isomorphic to the image of $\tilde\rho_i$. This completes the decomposition. The fact that there are not more irreducible representations comes from the fact that, from Proposition \firstdecomposition, any irreducible representation $V$ of $G$ is isomorphic to a subrepresentation of $K(G)$, and the uniqueness of the decomposition implies that $V$ must be one of the $V_i$'s.
\qed

\bigskip
\noindent{\bf \S\field. \tfield.} 
\medskip

In Theorem \ArtinWedderburn, we could expect the subalgebras $W_i$ to be isomorphic to $End(V_i)$, as Theorem \irreduciblesurjective\ asserts to happen when $K$ is algebraically closed. However, this is not the case in general. We will devote this section to discuss when this is so. Moreover, we will find an explicit method to find the representations $V_i$ from the idempotents, although the choice will not be canonical.

We start with a characterization of the surjectivity of $\tilde\rho$, of a different flavor of Theorem \irreduciblesurjective, since there we started with a representation $V$, while our goal now is precisely, given a centrally primitive idempotent $\beta$, to recover from $\beta K(G)$ the corresponding irreducible representation $V$. The idea will be to work with dimensions.

\proclaim Proposition \findrepresentation. Let $\rho:G\to Aut(V)$ be an irreducible representation of $G$ corresponding to an idempotent $\beta\in K(G)$ as in Theorem \ArtinWedderburn\ (in particular we are assuming $|G|$ is not divisible by the characteristic of $K$). Then the map $\tilde\rho:K(G)\to End(V)$ is surjective if and only if there exists a nonzero $\alpha\in\beta K(G)$ such that the right multiplication $\cdot\alpha:\beta K(G)\to\beta K(G)$ has rank at most $\sqrt{\dim\beta K(G)}$. Moreover, in this case, the representation $V$ is isomorphic to $\beta K(G)\alpha$.

\noindent{\it Proof.} Assume first that $\tilde\rho$ is surjective, so that $\beta K(G)$ is isomorphic to $End(V)$. Fix an endomorphism $f:V\to V$ of rank one. Hence the right composition $\circ f:End(V)\to End(V)$ has as image $V':=\{g\circ f\ |\ g\in End(V)\}\subset End(V)$. If $H\subset V$ is the kernel of $f$, then $V'$ is the set of endomorphisms of $V$ whose kernel contains $H$ and hence the representation $V'$ is isomorphic to $V$ (see Remark \decompositionnotunique). Hence the rank of $\circ f$ is the dimension of $V$, which is exactly the square root of the dimension of $End(V)$. 

Reciprocally, assume that we have a nonzero $\alpha\in\beta K(G)$ such that $\cdot\alpha$ has rank at most $\sqrt{\dim\beta K(G))}$. This implies that $\beta K(G)\alpha$ is a nonzero subrepresentation of $\beta K(G)$ of dimension at most $\sqrt{\dim\beta K(G))}$. Since $\beta K(G)$ is isomorphic to a subrepresentation of $End(V)$, it has dimension at most the square of $\dim(V)$, and it decomposes as direct sum of copies of $V$. This implies that $\beta K(G)\alpha$ is isomorphic to $V$, as wanted.
\qed

Let us apply the above criterion to a couple of examples.

\noindent{\bf Example \SchurSthree.} We have seen in Example \exampleSthree\ that $K(S_3)$ has the following centrally primitive idempotents
$$\beta_1={1+\alpha_{(2,1)}+\alpha_{(3)}\over6}={1+e_{(1\ 2)}+e_{(1\ 3)}+e_{(2\ 3)}+e_{(1\ 2\ 3)}+e_{(1\ 3\ 2)}\over6}$$
$$\beta_2={1-\alpha_{(2,1)}+\alpha_{(3)}\over6}={1-e_{(1\ 2)}-e_{(1\ 3)}-e_{(2\ 3)}+e_{(1\ 2\ 3)}+e_{(1\ 3\ 2)}\over6}$$
$$\beta_3={4-2\alpha_{(3)}\over6}={4-2e_{(1\ 2\ 3)}-2e_{(1\ 3\ 2)}\over6}.$$

It is clear that $\alpha\beta_1=\beta_1=\beta_1\alpha$ for any $\alpha\in K(S_3)$, so that $V_1:=\beta_1K(S_3)$ is the one-dimensional subspace spanned by $\beta_1$ and any $\alpha\in K(S_3)$ acts as the identity on $V_1$, so that $\beta_1$ induces the trivial representation.

Similarly, for any $\sigma\in S_3$, one has now $e_\sigma\beta_2=sgn(\sigma)\beta_2=\beta_2e_\sigma$. Hence we have again that $V_2:=\beta_2K(S_3)$ is the one dimensional subspace spanned by $\beta_2$, but now any $\sigma\in S_3$ acts as the multiplication by the sign of $\sigma$. Hence $\beta_2$ corresponds to the alternating representation.

We have to study finally the representation induced by $\beta_3$. By exclusion, it should be the two-dimensional standard representation (Example \standardirreducible), but we will try to construct it directly by hand from $\beta_3$. The first observation is that $V_3:=\beta_3K(S_3)$ has dimension four, hence it must correspond to $End(V)$ for some two-dimensional vector space $V$ that we want to find. The first problem is that there is no natural basis for $V_3$, so that we need to fix one. Instead of getting one multiplying $\beta_3$ by suitable elements $e_\sigma$, it is more convenient to find first implicit equations of $V_3$ inside $K(S_3)$. Observe that $V_3$ can be characterized as those elements of $K(S_3)$ whose product by $\beta_1$ and $\beta_2$ is zero. This is equivalent, adding and substracting, that the product with $1+e_{(1\ 2\ 3)}+e_{(1\ 3\ 2)}$ and $e_{(1\ 2)}+e_{(1\ 3)}+e_{(2\ 3)}$ is zero. Since this second element is obtained from the first one multiplying by $e_{(1\ 2)}$, we obtain that an element
$$\gamma=\mu_11+\mu_2e_{(1\ 2)}+\mu_3e_{(1\ 3)}+\mu_4e_{(2\ 3)}+\mu_5e_{(1\ 2\ 3)}+\mu_6e_{(1\ 3\ 2)}\in K(G)$$
is in $V_3$ if and only if
$$\mu_1+\mu_5+\mu_6=0$$
$$\mu_2+\mu_3+\mu_4=0.$$
We can thus take as a basis of $V_3$
$$v_1:=1-e_{(1 2 3)},\ v_2:=1-e_{(1 3 2)},\ v_3=e_{(1 2)}-e_{(1 3)},\ v_4:=e_{(1 2)}-e_{(2 3)}$$
and the multiplication is given by 
$$(y_1v_1+y_2v_2+y_3v_3+y_4v_4)(x_1v_1+x_2v_2+x_3v_3+x_4v_4)=z_1v_1+z_2v_2+z_3v_3+z_4v_4$$
where
$$z_1=(2x_1+x_2)y_1+(x_1-x_2)y_2+(x_3+2x_4)y_3+(-x_3+x_4)y_4$$
$$z_2=(-x_1+x_2)y_1+(x_1+2x_2)y_2+(x_3-x_4)y_3+(2x_3+x_4)y_4$$
$$z_3=(2x_3+x_4)y_1+(x_3-x_4)y_2+(x_1+2x_2)y_3+(-x_1+x_2)y_4$$
$$z_4=(-x_3+x_4)y_1+(x_3+2x_4)y_2+(x_1-x_2)y_3+(2x_1+x_2)y_4.$$
Hence, the right multiplication by $x_1v_1+x_2v_2+x_3v_3+x_4v_4$ has a rank at most two if and only if the matrix
$$\pmatrix{2x_1+x_2&x_1-x_2&x_3+2x_4&-x_3+x_4\cr
-x_1+x_2&x_1+2x_2&x_3-x_4&2x_3+x_4\cr
2x_3+x_4&x_3-x_4&x_1+2x_2&-x_1+x_2\cr
-x_3+x_4&x_3+2x_4&x_1-x_2&2x_1+x_2}$$
has rank at most two. One can see that this is equivalent to the vanishing of $x_1^2+x_1x_2+x_2^2-x_3^2-x_3x_4-x_4^2$ (which is in fact the square root of the determinant of the above matrix). From Proposition \findrepresentation, the representation $V$ we are looking for is isomorphic to the product of an element of rank one in $\beta_3K(S_3)$ with $K(S_3)$. For example, $v_2+v_3=1+e_{(1 2)}-e_{(1 3)}-e_{(1 3 2)}$ has rank one. This is the particular {\it Young symmetrizer} (see \S6.1 of \FultonHarris) used to understand all the irreducible representations of the symmetric groups. In fact, this is the standard way to prove that $\tilde\rho$ is surjective for any representation of any symmetric group (see also \LuxPahlings\ Theorem 3.3.8).

The reader maybe noticed that the coefficient of $1$ in the denominator of each $\beta_1$ above coincided exactly with the dimension of $\beta_iK(S_3)$. Of course this is not a coincidence, and it is implicitely related to the classical theory of characters (the character table of $S_3$ in Example \exampleSthree, hidden in the coordinates of the points or in the coefficients of the nilpotent elements, i.e. in the columns and rows respectively of the matrices $A$ and $B$). We will study this relation more closely.  For the time being, a first step is the following:

\proclaim Proposition \charactersidempotents. Let $G$ be a finite group and let $K$ be a field whose characteristic is not a divisor of $|G|$. Then:
\item{(i)} If $\alpha=\sum_g\mu_ge_g$, then each $\mu_g$ is $1\over|G|$ times the trace of the endomorphism $e_{g^{-1}}\alpha\cdot:K(G)\to K(G)$. 
\item{(ii)} The element $\beta=\sum_{g\in G}{1\over|G|}e_g$ is a centrally primitive idempotent whose associated irreducible representation is the trivial representation.
\item{(iii)} If the centrally primitive idempotent $\beta=\sum_{\lambda\in P}b_\lambda\alpha_\lambda$ corresponds to the irreducible representation $\rho:G\to Aut(V)$, then the coefficient $b_1$ of $1$ is $b_1={\dim\tilde\rho(K(G))\over|G|}$. 
\item{(iv)} If in (iii) we assume that $\tilde\rho$ is surjective, $b_\lambda={trace(\lambda^{-1})\dim V\over|G|}$ where $trace(\lambda^{-1})$ is the trace of any $g^{-1}\cdot:V\to V$ for $g\in\lambda$.

\noindent{\it Proof.} To prove (i), observe first that, for each $h\in G$, the image by $e_{g^{-1}}\alpha\cdot$ of $e_h$ is $\sum_{g'}\mu_{g'}e_{g^{-1}g'h}$, so that the coefficient of $e_h$ is $\mu_g$. Therefore the matrix of $e_{g^{-1}}\alpha\cdot$ with respect to the basis given by the vectors $e_h$ has only $\mu_g$ in its diagonal, so that its trace is $|G|\mu_g$.

Part (ii) is a consequence of (i). However, it can be proved easily by hand. Indeed, it is clear that $e_g\beta=\beta=\beta e_g$ for each $g\in G$. In particular, $\beta$ is in the center of $K(G)$, and $\alpha\beta=\beta$ for all $\alpha\in K(G)$. Hence $\beta^2=\beta$, so that $\beta$ is an idempotent element. Since $\beta K(G)$ is the one-dimensional subspace $V$ spanned by $\beta$, it follows that $\beta$ is necessarily primitive and the corresponding  representation is $V$, which is trivial since the multiplication by any $\alpha$ restricts to the identity on $V$.

For (iii), consider the endomorphism $\beta\cdot:K(G)\to K(G)$ consisting of the left multiplication by $\beta$. By part (i), $b_1={trace(\beta\cdot)\over|G|}$. Using the decomposition of Theorem \ArtinWedderburn, the endomorphism $\beta\cdot$ is the identity in the summand $\beta K(G)$, and zero in the summands corresponding to the other centrally primitive idempotents, so that its trace is the dimension of $\beta K(G)$.

For part (iv), the coefficient $b_\lambda$ of $\beta$ is the coefficient of any $e_g$ in (i) for $g\in\lambda$. Hence, by (i), $trace(e_{g^{-1}}\beta\cdot)=|G|b_\lambda$. On the other hand, as in the proof of (iii), $e_{g^{-1}}\beta\cdot$ is the multiplication by $g^{-1}$ on $\beta K(G)$ and zero on the other summands, so that its trace is the same as the trace of  $e_{g^{-1}}\beta\cdot:\beta K(G)\to\beta K(G)$. Since $\tilde\rho$ is surjective, Lemma \groupalgebra\ implies that this is the trace of $g^{-1}:End(V)\to End(V)$, hence $\dim V$ times the trace of $g^{-1}\cdot:V\to V$. 
\qed

The good behavior of the symmetric group in Example \SchurSthree\ is quite particular. In general, the situation is more complicated. We include the following illuminating example, who was suggested to us by J\"urgen M\"uller:

\noindent{\bf Example \quaternions.} Let us try to repeat the construction of Example \SchurSthree\ for the {\it quaternion group} $G=Q_8=\{\pm1,\pm i,\pm j, \pm k\}$, with $i^2=j^2=k^2=ijk=-1$. In this case, the center of $K(Q_8)$ is generated by the elements
$$1=e_1$$
$$\alpha_1=e_{-1}$$
$$\alpha_2=e_i+e_{-i}$$
$$\alpha_3=e_j+e_{-j}$$
$$\alpha_4=e_k+e_{-k}$$
and one easily checks $ZK(Q_8)\cong K[x,y,z,w]/I$ where 
$$I=(x^2-1,xy-y,xz-z,xw-w,y^2-2x-2,yz-2w,yw-2z,z^2-2x-2,zw-2y,w^2-2x-2)$$
which is the ideal of
$$V(I)=\{(-1,0,0,0),(1,2,2,2),(1,2,-2,-2),(1,-2,2,-2),(1,-2,-2,2)\}.$$
The evaluation at these points provides an isomorphism $ZK(Q_8)\to K^5$ whose inverse yields
$$(1,0,0,0,0)\mapsto\beta_1:={4-4\alpha_1\over8}$$
$$(0,1,0,0,0)\mapsto\beta_2:={1+\alpha_1+\alpha_2+\alpha_3+\alpha_4\over8}$$
$$(0,0,1,0,0)\mapsto\beta_3:={1+\alpha_1+\alpha_2-\alpha_3-\alpha_4\over8}$$
$$(0,0,0,1,0)\mapsto\beta_4:={1+\alpha_1-\alpha_2+\alpha_3-\alpha_4\over8}$$
$$(0,0,0,0,1)\mapsto\beta_5:={1+\alpha_1-\alpha_2-\alpha_3+\alpha_4\over8}$$
As Proposition \charactersidempotents(iii) shows, the representations given by $\beta_2,\beta_3,\beta_4,\beta_5$ are one-dimensional, so there is nothing interesting to study, while $\beta_1$ is expected to provide an irreducible two-dimensional representation. Let us check directly what happens in this case. As one easily checks, a basis of $\beta_1K(Q_8)$ is 
$$v_1:=e_1-e_{-1},\ v_2:=e_i-e_{-i},\ v_3:=e_j-e_{-j},\ v_4:=e_k-e_{-k}$$
and the multiplication is given by 
$$(y_1v_1+y_2v_2+y_3v_3+y_4v_4)(x_1v_1+x_2v_2+x_3v_3+x_4v_4)=z_1v_1+z_2v_2+z_3v_3+z_4v_4$$
where
$$z_1=2x_1y_1-2x_2y_2-2x_3y_3-2x_4y_4$$
$$z_2=2x_2y_1+2x_1y_2+2x_4y_3-2x_3y_4$$
$$z_3=2x_3y_1-2x_4y_2+2x_1y_3+2x_2y_4$$
$$z_4=2x_4y_1+2x_3y_2-2x_2y_3+2x_1y_4.$$
Hence, the right multiplication by $x_1v_1+x_2v_2+x_3v_3+x_4v_4$ has a two-dimensional linear space in its kernel if and only if the matrix
$$\pmatrix{x_1&-x_2&-x_3&-x_4\cr
x_2&x_1&x_4&-x_3\cr
x_3&-x_4&x_1&x_2\cr
x_4&x_3&-x_2&x_1}$$
has rank at most two. However, this is equivalent to the vanishing of $x_1^2+x_2^2+x_3^2+x_4^2$ (the square root of the determinant of the above matrix), which is never zero if, for example, $K={\Bbb Q}$ (hence, when $K={\Bbb Q}$ the irreducible representation corresponding to $\beta_3$ is the whole $\beta_3K(Q_8)$, which has dimension $4$). If, instead, $K$ contains $\sqrt{-1}$, then a solution is $x_1=x_2=0$, $x_3=1$, $x_4=\sqrt{-1}$, since it provides a skew-symmetric matrix whose Pfaffian is zero. This is not the only possibility. For example, if $K$ contains a primitive cubic root $\omega$ of 1, then also $x_1=0$, $x_2=1$, $x_3=\omega$, $x_4=\omega^2$ would be another solution. Observe that, in characteristic $p$, since any prime number is a sum of four squares, it follows that $x_1^2+x_2^2+x_3^2+x_4^2=0$ always have a solution, so that we find an irreducible representation of dimension two.

\noindent{\bf Example \dihedral.} It is also interesting to compare the above example with the dihedral group $D_4$. Interpreting it as the group of               isometries of the square of consecutive vertices labelled as $1,2,3,4$, we can identify it with the subgroup of the symmetric group $G=D_4=\{(1),(1\ 3)(2\ 4),(1\ 3),(2\ 4)),(1\ 2\ 3\ 4),(1\ 4\ 3\ 2),(1\ 2)(3\ 4),(1\ 4)(2\ 3)\}$. Now the center of $K(D_4)$ is generated by the elements
$$1=e_{(1)}$$
$$\alpha_1=e_{(1\ 3)(2\ 4)}$$
$$\alpha_2=e_{(1\ 3)}+e_{(2\ 4)}$$
$$\alpha_3=e_{(1\ 2\ 3\ 4)}+e_{(1\ 4\ 3\ 2)}$$
$$\alpha_4=e_{(1\ 2)(3\ 4)}+e_{(1\ 4)(2\ 3)}$$
and they satisfy exactly the same relations as the corresponding elements of $K(Q_8)$. We thus get the same set of points $V(I)$ and the same idempotent elements. In particular, $\beta_1={4-4\alpha_1\over8}$ will be the only idempotent giving a representation of dimension bigger than one. Now, a basis of $\beta_1K(D_4)$ is 
$$v_1:=e_{(1)}-e_{(1\ 3)(2\ 4)},\ v_2:=e_{(1\ 3)}-e_{(2\ 4)},\ v_3:=e_{(1\ 2\ 3\ 4)}-e_{(1\ 4\ 3\ 2)},\ v_4:=e_{(1\ 2)(3\ 4)}-e_{(1\ 4)(2\ 3)}$$
and the multiplication is given by 
$$(y_1v_1+y_2v_2+y_3v_3+y_4v_4)(x_1v_1+x_2v_2+x_3v_3+x_4v_4)=z_1v_1+z_2v_2+z_3v_3+z_4v_4$$
where
$$z_1=2x_1y_1+2x_2y_2-2x_3y_3+2x_4y_4$$
$$z_2=2x_2y_1+2x_1y_2+2x_4y_3-2x_3y_4$$
$$z_3=2x_3y_1+2x_4y_2+2x_1y_3-2x_2y_4$$
$$z_4=2x_4y_1+2x_3y_2-2x_2y_3+2x_1y_4.$$
We obtain now that the right multiplication by $x_1v_1+x_2v_2+x_3v_3+x_4v_4$ has a two-dimensional linear space in its kernel if and only if the matrix
$$\pmatrix{x_1&x_2&-x_3&x_4\cr
x_2&x_1&x_4&-x_3\cr
x_3&x_4&x_1&-x_2\cr
x_4&x_3&-x_2&x_1}$$
has rank at most two. This is equivalent to the vanishing of $x_1^2-x_2^2+x_3^2-x_4^2$ (again the square root of the determinant of the above matrix), and we can find a solution (for example $v_1+v_2$) without need of extending the field. As a consequence, even if the set of points is the same for $Q_8$ and $D_4$, the set of representations of both groups is not the ``same'' for an arbitrary field, but only in the case of a splitting field.

We summarize here and make more explicit what is our main construction in this section, fixing the notation that we will use in the rest of the note. 

\noindent{\bf Remark \mainconstruction.} Given a finite group $G$ and fixing a characteristic that is not a divisor of $|G|$, we consider the minimal field $K_0$ with such characteristic (i.e. $K_0={\Bbb Z}_p$ or $K_0={\Bbb Q})$. Identify $ZK_0(G)$ with the quotient $K_0[\{x_\lambda\}_{\lambda\in P'}]/I_0$ of Remark \idealcenter. Then any field $K$ containing the finite extension of $K_0$ attaching the coordinates of the points of $V(I_0)$ works fine for the main decomposition in Theorem \ArtinWedderburn. In particular, there is a bijection between the set $Irr_K(G)$ of irreducible representations of $G$ and the set of centrally primitive idempotents of $K(G)$; we will write $\beta_V$ for the idempotent corresponding to $V\in Irr_K(G)$). On the other hand, using the notation of Example \exampleSthree, we proved that we have an isomorphism 
$$\varphi:ZK(G)\cong K[\{x_\lambda\}_{\lambda\in P'}]/I\to K^{|P|}$$ 
consisting of the evaluation at the set of $|P|$ points $V(I)\subset\A^{|P|-1}$. Since each centrally primitive idempotent $\beta_V$ maps by $\varphi$ to an element of the canonical basis of $K^{|P|}$, there is only one point $a_V\in V(I)$ such that that is not zero when evaluated at $\beta_V$. Hence we will write $V(I)=\{a_V\}_{V\in Irr_K(G)}$ and their coordinates will be $a_{V,\lambda}$, when $\lambda\in P$. Sometimes we will add a first coordinate $a_{V,1}=1$ corresponding to the conjugate class of $1$. Using the basis $\{\alpha_\lambda\}_{\lambda\in P}$ in $ZK(G)$,  the isomorphism $\varphi$ is given by 
$$\sum_{\lambda\in P}\mu_\lambda\alpha_\lambda=\big(\{\mu_\lambda\}_{\lambda\in P}\big)\ \big(\{\alpha_\lambda\}_{\lambda\in P}\big)^t\mapsto\big(\{\mu_\lambda\}_{\lambda\in P}\big)A$$
where $A$ is the matrix whose columns are the coordinates $a_{V,\lambda}$ (including $a_{V,1}=1$) of the points of $V(I)$. The inverse of $\varphi$ is thus given by 
$$\big(\{\nu_\lambda\}_{\lambda\in P}\big)\mapsto\big(\{\nu_\lambda\}_{\lambda\in P}\big)\ B\ \big(\{\alpha_\lambda\}_{\lambda\in P}\big)^t$$
where $B:=A^{-1}$ is the matrix whose rows are the coefficients of the centrally primitive idempotents $\{\beta_V\}_{V\in Irr_K(G)}$ of $K(G)$. We stress that, although there is the same number of irreducible representations of $G$ and of conjugacy classes of $G$, there is not an a priori intrinsic bijection (when $G$ is a symmetric group, an explicit bijection is given via the Young symmetrizers; see \FultonHarris\ Theorem 4.3 for details). 

The remaining point is whether we have the surjectivity for the maps $\tilde\rho:K(G)\to End(V)$ corresponding to irreducible representations. Using Proposition \findrepresentation, and as we have seen in the previous examples, this depends on an algebraic set to have points with coordinates in $K$. Since this is true for the algebraic closure of $K$ (Theorem \irreduciblesurjective), we know that we can find such points with coordinates in a finite extension of $K$. Example \quaternions\ shows that we have several choices for a minimal such extension, since in general we can choose our points in different extensions of $K$. The general definition is:

\noindent{\bf Definition.} A {\it splitting field} of a group is a field $K$ such that, for each irreducible representation $\rho:G\to Aut(V)$, the corresponding map $\tilde\rho:K(G)\to End(V)$ is surjective.

When we are working with a splitting field, we have now the good properties:

\proclaim Theorem \propertiessplitting. Let $G$ be a finite group and let $K$ be a splitting field of $G$ whose characteristic does not divide $|G|$. Then:
\item{(i)} {\rm (Artin-Wedderburn decomposition)} $K(G)\cong\bigoplus_iEnd(V_i)$, where the $V_i$'s are the irreducible representations of $G$.
\item{(ii)} {\rm (Schur)} Any endomorphism of an irreducible representation that is a morphism of representations is necessarily the multiplication by a constant.
\item{(iii)} If $V$ is an irreducible representation, for any $\alpha\in ZK(G)$, let $\lambda_\alpha\in K$ be the constant such that, according to (ii), the morphism $\alpha\cdot:V\to V$ is the multiplication by $\lambda_\alpha$. Then the map $\omega_V:ZK(G)\to K$ defined by $\alpha\mapsto\lambda_\alpha$ is a homomorphism of $K$-algebras. 
\item{(iv)} If $K\subset K'$ is any field extension and $V$ is a representation of $G$ over $K$, let $V'=V\otimes_KK'$ be the extended representation over $K'$. Then $V$ is irreducible if and only if $V'$ is irreducible.

\noindent{\it Proof.} Part (i) is just Theorem \ArtinWedderburn\ using that $K$ is a splitting field.

For part (ii), let $\varphi:V\to V$ be a morphism of representations, with $V$ irreducible. This means that $\varphi$ commutes with any endomorphism of $V$  in the image of $\rho$, hence with any endomorphism in the image of $\tilde\rho$. Since $K$ is a splitting field, this means that $\varphi$ commutes with any endomorphism of $V$, so that it is necessarily the multiplication by a constant.

For part (iii), consider $\alpha_1,\alpha_2\in ZK(G)$, so that $\alpha_1v=\lambda_{\alpha_1}v$ and $\alpha_2v=\lambda_{\alpha_2}v$ for any $v\in V$. Thus 
$$(\alpha_1+\alpha_2)v=(\lambda_{\alpha_1}+\lambda_{\alpha_2})v$$
so that $\lambda_{\alpha_1+\alpha_2}=\lambda_{\alpha_1}+\lambda_{\alpha_2}$, i.e. $\omega_V(\alpha_1+\alpha_2)=\omega_V(\alpha_1)+\omega_V(\alpha_2)$;  and also
$$(\alpha_1\alpha_2)v=\alpha_1(\lambda_{\alpha_2}v)=\lambda_{\alpha_1}(\lambda_{\alpha_2}v)$$
so that $\lambda_{\alpha_1\alpha_2}=\lambda_{\alpha_1}\lambda_{\alpha_2}$, i.e. $\omega_V(\alpha_1\alpha_2)=\omega_V(\alpha_1)\omega_V(\alpha_2)$.

Finally, for part (iv), it is clear that the irreducibility of $V'$ implies the irreducibility of $V$, since any subrepresentation of $V$ would extend to $V'$. Reciprocally, if $V$ is irreducible, we know that $\tilde\rho:K(G)\to End(V)$ is surjective. Tensoring with $K'$, we see that the corresponding $\tilde\rho':K'(G)\to End(V')$ is also surjective, so that $V'$ is irreducible.
\qed

With these results, we can describe the coordinates of the point corresponding to an irreducible representation:

\proclaim Proposition \characterspoints. Let $G$ be a finite group and let $K$ be a splitting field whose characteristic is not a divisor of $|G|$. Let $V$ be an irreducible representation and let $\omega_V:ZK(G)\to K$ be the homomorphism defined in Theorem \propertiessplitting(iii). Then:
\item{(i)} With the identification $ZK(G)\cong K[\{x_\lambda\}_{\lambda\in P'}]/I$ of Remark \idealcenter, the map $\omega_V$ is the evaluation of classes of polynomials at the point $a_V$.
\item{(ii)} The $\lambda$-coordinate of $a_V$ is $a_{V,\lambda}={trace(\lambda)|\lambda|\over\dim(V)}$, where $trace(\lambda)$ is the trace of the multiplication map $g\cdot:V\to V$ for any $g\in\lambda$.

\noindent{\it Proof.} For part (i), let $\beta_V$ be the centrally primitive idempotent corresponding to $V$. Then the multiplication by $\beta_V$ is the identity on $V$, so that $\omega_V(\beta_V)=1$. If instead $V'$ is a different irreducible representation and $\beta_{V'}$ is its corresponding idempotent, then the multiplication by $\beta_{V'}$ is the zero map on $V$, so that $\omega_V(\beta_{V'})=0$. This proves that, if we consider the homomorphism of $K$-algebras $ZK(G)\to\Pi_{V'}K$ given by the product of the maps $\omega_{V'}$, the image of each $\beta_{V'}$ is the element of the canonical basis whose nonzero coordinate is the one corresponding to $V'$. Hence this homomorphism is the isomorphism given in Proposition \evaluationisomorphism(iv) consisting of the evaluation at the points of $V(I)$. In particular, $\omega_V$ is the evaluation at $a_V$, as wanted.

For part (ii), observe that part (i) implies that the $\lambda$-coordinate of $a_V$ is $\omega_V(\alpha_\lambda)$. This means that the map $\alpha_\lambda\cdot:V\to V$ is the multiplication by $a_{V,\lambda}$, so that its trace is $\dim V$ times $a_{V,\lambda}$. On the other hand, since  $\alpha_\lambda=\sum_{g\in\lambda}e_g$, the trace of $\alpha_\lambda\cdot$ is the sum of the traces of the multiplications by all $g\in\lambda$, hence equal to $trace(\lambda)|\lambda|$ (all the traces of the multiplications by $g\in\lambda$ are the same because they are conjugate). Hence the result follows.
\qed

We can now recover the main result in the classical theory of representations. We start we a definition:

\noindent{\bf Definition.} The {\it character of a representation} $V$ of a group $G$ is the map $\chi_V:G\to K$ mapping each $g$ to the trace of the map $g\cdot:V\to V$. As already remarked, we could define  $\chi_V:P\to K$, and write $\chi_V(\lambda)$ for $\chi_V(g)$ for any $g\in\lambda$.

\proclaim Corollary \duality. Given two irreducible representations $V,W$ of a finite group $G$,
$${1\over|G|}\sum_{g\in G}\chi_V(g^{-1})\chi_W(g)=\left\{\matrix{1\hbox{ if } V=W\cr 0\hbox{ if } V\ne W.}\right.$$

\noindent{\it Proof.} Regrouping addends according to their conjugation classes, the above sum becomes
$${1\over|G|}\sum_{g\in G}\chi_V(g^{-1})\chi_W(g)={1\over|G|}\sum_{\lambda\in P}|\lambda|\chi_V(\lambda^{-1})\chi_W(\lambda)=$$
$$={\dim W\over\dim V}\sum_{\lambda\in P}{\chi_V(\lambda^{-1})\dim V\over|G|}{\chi_W(\lambda)|\lambda|\over\dim W}=\left\{\matrix{1\hbox{ if } V=W\cr 0\hbox{ if } V\ne W}\right.$$
where the last equality comes from the fact that we are multiplying the $V$-th row of the matrix $B$ corresponding to the coordinates of the idempotent $\beta_V$ (see Proposition \charactersidempotents) and the $W$-the column of its inverse $A$ corresponding to the coordinates of the point $a_W$ (see Proposition \characterspoints).
\qed

We finish this section by observing that we have also a way of decomposing representations, and we can do it using only the coordinates of the points $a_V$ (we will use this result in the next section for an arbitrary vector space $W$ on which $G$ acts, not necessarily of finite dimension):

\proclaim Proposition \basisevaluation. Let $W$ be any representation of $G$. Then $W=\bigoplus_{V\in Irr_K(G)}\beta_VW$, and each $\beta_VW$ is the direct sum of the irreducible representation $V$ a finite number of times. Moreover, the subrepresentation $\beta_VW$ can be characterized as the set of vectors $w\in W$ satisfying one of the following conditions:
\item{(i)} $\alpha_\lambda w=a_{V,\lambda}w$ for each conjugation class $\lambda$ different from the one of $1$.
\item{(ii)} $\alpha w=0$, where $\alpha$ is a fixed element of $ZK(G)$ of the form $\alpha=\sum_{V'\ne V}\mu_{V'}\beta_{V'}$ with all $\mu_{V'}\ne0$.

\noindent{\it Proof.} The decomposition $W=\bigoplus_V\beta_VW$ is immediate. Since $\beta_{V'}(\beta_VW)=0$ if $V'\ne V$, it follows that the only irreducible subrepresentation that $\beta_VW$ can contain is precisely $V$. Moreover, $\beta_VW$ is characterized as the set of vectors $w\in W$ such that $\beta_{V'}w=0$ for any irreducible representation $V'\not\cong V$. 

For characterization (i), observe that all these $\beta_{V'}$ form a basis of the kernel of the map $ZK(G)\to K$ that corresponds to the evaluation map $K[\{x_\lambda\}_{\lambda\in P'}]/I\to K$ at the point $a_V$. Since the elements $\alpha_\lambda-a_{V,\lambda}e_1$ form another basis of the kernel, characterization (i) follows.

For characterization (ii), it is obvious that, if $\beta_{V'}w=0$ for any irreducible representation $V'\not\cong V$, then also $\alpha w=0$. Reciprocally, if $\alpha w=0$, we get that, for any irreducible representation $V'\not\cong V$,
$$0=\beta_{V'}\alpha w=\beta_{V'}\sum_{V''\ne V}\mu_{V''}\beta_{V''}w=\mu_{V'}\beta_{V'}w$$
so that, since $\mu_{V'}\ne0$, we get $\beta_{V'}w=0$.
\qed

\noindent{\bf Remark \remarkequations.} Characterization (ii) could look difficult to use, specially if we just want to use the coordinates of the points $a_V$.  However, this is not the case, since, with the identifications $ZK(G)\cong K[\{x_\lambda\}_{\lambda\in P'}]/I\cong K^{|P|}$, giving an element $\alpha$ as in Proposition \basisevaluation\ is equivalent to giving the equation $f\in K[\{x_\lambda\}_{\lambda\in P'}]$ passing through $a_V$ but not through any other point $a_{V'}$. For instance, in Example \exampleSthree\ we had the points $(3,2)$, $(-3,2)$ and $(0,-1)$. We can thus take polynomials $x-3$, $x+3$ and $y+1$ to separate the three points, hence improving characterization (i). These correspond respectively to elements (still with the notation of  Example \exampleSthree)
$$\alpha_{(2,1)}-3=e_{(1\ 2)}+e_{(1\ 3)}+e_{(2\ 3)}-3e_{(1)}=(e_{(1\ 2)}-e_{(1)})+(e_{(1\ 3)}-e_{(1)})+(e_{(2\ 3)}-e_{(1)})$$
$$\alpha_{(2,1)}+3=e_{(1\ 2)}+e_{(1\ 3)}+e_{(2\ 3)}+3e_{(1)}=(e_{(1\ 2)}+e_{(1)})+(e_{(1\ 3)}+e_{(1)})+(e_{(2\ 3)}+e_{(1)})$$
$$\alpha_{(3)}=e_{(1)}+e_{(1\ 2 \ 3)}+e_{(1\ 3\ 2)}.$$
Observe that the first (resp. second) element corresponds to the sum of the three standard conditions for a symmetric (resp. skew-symmetric) function. In fact, the elements $w$ of subspace $\beta_1W$ (resp. $\beta_2W$) can also be characterized by $e_{(1\ 2)}w=e_{(1\ 3)}w=e_{(2\ 3)}w=w$ (resp. $e_{(1\ 2)}w=e_{(1\ 3)}w=e_{(2\ 3)}w=-w$), because $\beta_1$ (resp. $\beta_2$) is killed by $e_{(1\ 2)}-e_{(1)}$, $e_{(1\ 3)}-e_{(1)}$, $e_{(2\ 3)}-e_{(1)}$ (resp. $e_{(1\ 2)}+e_{(1)}$, $e_{(1\ 3)}+e_{(1)}$, $e_{(2\ 3)}+e_{(1)}$).

\bigskip
\noindent{\bf \S\symmetries. \tsymmetries.} 
\medskip

We want to apply what we have done to understand symmetries of functions in $d$ variables. Our approach to study the irreducible representations from the points associated to the group will allow to find easy equations for the possible symmetries. We start explaining why both things are related.

\noindent{\bf Example \Jsymmetry.} Everyone knows how to decompose a function in two variables into a sum of its symmetric part and its skew-symmetric part. In three variables, a function $F:X\times X\times X\to K$ has a symmetric part
$$F_s(x_1,x_2,x_3):=\matrix{{F(x_1,x_2,x_3)+F(x_1,x_3,x_2)+F(x_2,x_1,x_3)+F(x_2,x_3,x_1)+F(x_3,x_1,x_2)+F(x_3,x_2,x_1)\over6}},$$
a skew-symmetric part
$$F_a(x_1,x_2,x_3):=\matrix{{F(x_1,x_2,x_3)-F(x_1,x_3,x_2)-F(x_2,x_1,x_3)+F(x_2,x_3,x_1)+F(x_3,x_1,x_2)-F(x_3,x_2,x_1)\over6}}$$
and, to complete the decomposition we need another summand
$$F_r(x_1,x_2,x_3):=\matrix{{4F(x_1,x_2,x_3)-2F(x_2,x_3,x_1)-2F(x_3,x_1,x_2)\over6}.}$$
The reader will immediately recognize the coefficients in the decomposition as the coefficients of the centrally primitive idempotents $\beta_1,\beta_2,\beta_3$ found in Example \SchurSthree. Moreover, a function as $F_r$ satisfies the condition  
$$F(x_1,x_2,x_3)+F(x_2,x_3,x_1)+F(x_3,x_1,x_2)=0$$
(this condition is sometimes called {\it $J$-symmetry}, because of the Jacobi identity of Lie brackets). This decomposition coincides with what we did in Remark \remarkequations. To formalize all this, consider the following:

\noindent{\bf Set-up.} Let $X$ be any set, $d$ a positive integer and let $K$ be an arbitrary field of characteristic zero (although it could be assumed to be positive characteristic if it is bigger than $d$). The symmetric group $S_d$ acts naturally on the set of  functions of $d$ variables $X\times\dots\times X\to K$ by defining, for any $\sigma\in S_d$ and any function $F$, the function $\sigma F$ as
$$(\sigma F)(x_1,\dots,x_d):=F(x_{\sigma(1)},\dots,x_{\sigma(d)}).$$ The space of these functions has a natural structure of vector space over $K$, and we will fix a linear subspace $C^d(X,K)$ invariant under the action of $S_d$. There is a natural way of extending the action of $S_d$ on $C^d(X,K)$ to an action of the elements of $K(S_d)$, by defining, for any $\alpha=\sum_{\sigma\in S_d}a_\sigma e_\sigma\in K(S_d)$ and $F\in C^d(X,K)$, a new function $\alpha F\in C^d(X,K)$ by
$$(\alpha F)(x_1,\dots,x_d):=\sum_{\sigma\in S_d}a_\sigma F(x_{\sigma(1)},\dots,x_{\sigma(d)}).$$

Let us put in our language some representation theory of the symmetric group:

\proclaim Proposition \idempotentssymmetric. If $G=S_d$ is the group of permutations of $d$ elements and $K$ is a field whose characteristic is not a divisor of $d!$:
\item{(i)} The element $\beta={1\over d!}\sum_{\sigma\in S_d}e_\sigma$ is the centrally primitive idempotent corresponding to the trivial representation.
\item{(ii)} The element $\beta={1\over d!}\sum_{\sigma\in S_d}sgn(\sigma)e_\sigma$ is the centrally primitive idempotent corresponding to the alternating representation.
\item{(iii)} The element $\beta={d-1\over d!}\sum_{\sigma\in S_d}(|Fix(\sigma)|-1)e_\sigma$, where $Fix(\sigma)$ is the set of elements fixed by $\sigma$, is the centrally primitive idempotent corresponding to the standard representation.

\noindent{\it Proof.} Part (i) is part (ii) of Proposition \charactersidempotents, and part (ii) comes from part (iv) of Proposition \charactersidempotents.

For part (iii), recall from Example \standardirreducible\ that the standard representation $V$ appears in the decomposition $K^d=V\oplus V'$, where $S_d$ acts on $K^d$ by permuting in the natural way the vectors of the canonical basis $e_1,\dots,e_d$, $V$ is the hyperplane of equation $x_1+\dots+x_d=0$ and $V'$ is the line spanned by $e_1+\dots+e_d$. Observe that the trace of the action of a permutation $\sigma\in S_d$ on $K^d$ is $|Fix(\sigma)|$ (which is also the number of points fixed by $\sigma^{-1}$, and it is also equal to the number of cycles of length one in the decomposition of $\sigma$ into disjoint cycles). Since the action of any $\sigma$ on $V'$, its trace is one, it follows that the trace of the multiplication by $\sigma^{-1}$ on $V$ is $|Fix(\sigma)|-1$. Since $V$ has dimension $d-1$, we conclude now from Proposition \charactersidempotents(iv).
\qed

Recalling our partition notation for the symmetric group, we can complete Proposition \idempotentssymmetric:

\proclaim Proposition \coordinatessymmetric. If $G=S_d$ is the group of permutations of $d$ elements and $K$ is a field whose characteristic is not a divisor of $d!$:
\item{(i)} The point corresponding to the trivial representation has as $\lambda$-coordinate $|\lambda|$.
\item{(ii)} The $\lambda$-coordinate of the point corresponding to the alternating representation is $sgn(\lambda)|\lambda|$, where $sgn(\lambda)$ is the sign of any permutation with partition $\lambda$.
\item{(iii)} The $\lambda$-coordinate of the point associated with the standard representation is \ $|\lambda|\big(|Fix(\lambda)|-1\big)\over d-1$, where $|Fix(\lambda)|$ is the number of $1$'s in the partition $\lambda$.

\noindent{\it Proof.} It is a consequence of Proposition \characterspoints. For part (iii) one also needs to recall from the proof of Proposition \idempotentssymmetric(iii) that the trace of the multiplication by $\sigma$ is $|Fix(\sigma)|-1$ and one easily observes that $|Fix(\sigma)|$ is the number of cycles of length one in its decomposition into disjoint cycles.
\qed 

Let us study in full detail the case of four variables:

\noindent{\bf Example \quadrilinearrepresentations.} If $d=4$, in $S_4$ we have now five conjugacy classes of permutations corresponding to the five possible partitions on $4$. More precisely, $ZK(S_4)$ is generated by
$$1:=e_{(1)}=\alpha_{(1,1,1,1,1)},\ \ \ \ \alpha_{(2,1,1)},\ \ \ \ \alpha_{(3,1)},\ \ \ \ \alpha_{(4)},\ \ \ \ \alpha_{(2,2)}.$$
Computing as in Example \exampleSthree\ the relations about the generators, we get that $ZK(S_4)$ is the quotient of $K[x,y,z,w]$ by the ideal $I$ generated by $x^2-6-3y-2w, xy-4x-4z,y^2-8-4y-8w,xz-3y-4w,yz-4x-4z,z^2-6-3y-2w,xw-x-2z,yw-3y,zw-2x-z,w^2-3-2w$. Now $I$ is the ideal of the points 
$$(x,y,z,w)=(6,8,6,3),(-6,8,-6,3),(0,-4,0,3),(-2,0,2,-1),(2,0,-2,-1)$$
so that we get an isomorphism of algebras of $ZK(S_4)$ with $K^5$ by evaluating at those points given by the matrix
$$A:=\pmatrix{1&6&8&6&3\cr1&-6&8&-6&3\cr1&0&-4&0&3\cr1&-2&0&2&-1\cr1&2&0&-2&-1}.$$
Hence the inverse is defined by
$$B=(A^{-1})^t=\pmatrix{1\over24&1\over24&1\over24&1\over24&1\over24\cr\cr
1\over24&-1\over24&1\over24&-1\over24&1\over24\cr\cr
4\over24&0&-2\over24&0&4\over24\cr\cr
9\over24&-3\over24&0&3\over24&-3\over24\cr\cr
9\over24&3\over24&0&-3\over24&-3\over24}$$
so that the centrally primitive idempotents are
$$\beta_s:={1+\alpha_{(2,1,1)}+\alpha_{(3,1)}+\alpha_{(4)}+\alpha_{(2,2)}\over24}$$
$$\beta_a:={1-\alpha_{(2,1,1)}+\alpha_{(3,1)}-\alpha_{(4)}+\alpha_{(2,2)}\over24}$$
$$\beta_1:={4-2\alpha_{(3,1)}+4\alpha_{(2,2)}\over24}$$
$$\beta_2:={9-3\alpha_{(2,1,1)}+3\alpha_{(4)}-3\alpha_{(2,2)}\over24}$$
$$\beta_3:={9+3\alpha_{(2,1,1)}-3\alpha_{(4)}-3\alpha_{(2,2)}\over24}$$
(according to Proposition \idempotentssymmetric(iii), this last idempotent corresponds to the standard representation). Applying this to $C^4(X,K)$, this corresponds to a decomposition $F=F_s+F_a+F_1+F_2+F_3$, where
$$F_s(x_1,x_2,x_3,x_4):={\sum_{\sigma\in S_d}F(x_{\sigma(1)},x_{\sigma(2)},x_{\sigma(3)},x_{\sigma(4)})\over24}$$
$$F_a(x_1,x_2,x_3,x_4):={\sum_{\sigma\in S_d}{\rm sgn}(\sigma)F(x_{\sigma(1)},x_{\sigma(2)},x_{\sigma(3)},x_{\sigma(4)})\over24}$$
$$\matrix{F_1(x_1,x_2,x_3,x_4):=\hfill\cr\cr
={2F(x_1,x_2,x_3,x_4)+2\sum_{\sigma\in(2,2)}F(x_{\sigma(1)},x_{\sigma(2)},x_{\sigma(3)},x_{\sigma(4)})-\sum_{\sigma\in(3,1)}F(x_{\sigma(1)},x_{\sigma(2)},x_{\sigma(3)},x_{\sigma(4)})\over12}}$$
$$\matrix{F_2(x_1,x_2,x_3,x_4):=\hfill\cr\cr$$
$$={3F(x_1,x_2,x_3,x_4)-\sum_{\sigma\in(2,1,1)\cup(2,2)}F(x_{\sigma(1)},x_{\sigma(2)},x_{\sigma(3)},x_{\sigma(4)})+\sum_{\sigma\in(4)}F(x_{\sigma(1)},x_{\sigma(2)},x_{\sigma(3)},x_{\sigma(4)})\over8}}$$
$$\matrix{F_3(x_1,x_2,x_3,x_4):=\hfill\cr\cr$$
$$={3F(x_1,x_2,x_3,x_4)-\sum_{\sigma\in(4)\cup(2,2)}F(x_{\sigma(1)},x_{\sigma(2)},x_{\sigma(3)},x_{\sigma(4)})+\sum_{\sigma\in(2,1,1)}F(x_{\sigma(1)},x_{\sigma(2)},x_{\sigma(3)},x_{\sigma(4)})\over8}.}$$
Clearly, $F_s$ is the symmetric part and $F_a$ is the alternating part (this can be done as in Remark \remarkequations). We study separately the rest of the components (again as in Remark \remarkequations) and compare with \MRS:

1) For the equations for the space of functions of type $F=F_1$ are given, we observe that we can take the hyperplane $y+4$ to characterize the point $(0,-4-0-3)$ associated to this symmetry. This corresponds to the element $\alpha_{(3,1)}+4e_{(1)}$, which we can write as the sum of the elements 
$$e_{(1)}+e_{(1\ 2\ 3)}+e_{(1\ 3\ 2)}$$
$$e_{(1)}+e_{(1\ 2\ 4)}+e_{(1\ 4\ 2)}$$
$$e_{(1)}+e_{(1\ 3\ 4)}+e_{(1\ 4\ 3)}$$
$$e_{(1)}+e_{(2\ 3\ 4)}+e_{(2\ 4\ 3)}.$$
Since the product of each of these elements with $\beta_1$ is zero, we conclude that we can characterize this kind of symmetry by the conditions
$$F(x_1,x_2,x_3,x_4)+F(x_2,x_3,x_1,x_4)+F(x_3,x_1,x_2,x_4)=0$$
$$F(x_1,x_2,x_3,x_4)+F(x_2,x_4,x_3,x_1)+F(x_4,x_3,x_1,x_2)=0$$
$$F(x_1,x_2,x_3,x_4)+F(x_3,x_2,x_4,x_1)+F(x_4,x_2,x_1,x_3)=0$$
$$F(x_1,x_2,x_3,x_4)+F(x_1,x_4,x_2,x_3)+F(x_1,x_3,x_4,x_2)=0$$
which are equations obtained in \MRS\ for case 3 in pages 851-852 (although the sum of them would be enough for the characterization).

2) We study now the space of functions of type $F=F_2$. Although one simple equation would be enough (for example, using the hyperplane $x+2=0$), we will try to find the characterization of \MRS. For this, we characterize the point $(-2,0,2,-1)$ using not only the hyperplane $x+w+3=0$ (because it contains also the point $(-6,8,-6,3)$ corresponding to the alternating part) but also the hyperplane $-x+y-z+w+1=0$. The element $\alpha_{(2,1,1)}+\alpha_{(2,2)}+3e_{(1)}$ corresponding to the first hyperplane can be written as the sum of the elements
$$e_{(1)}+e_{(1\ 2)}+e_{(3\ 4)}+e_{(1\ 2)(3\ 4)}$$
$$e_{(1)}+e_{(1\ 3)}+e_{(2\ 4)}+e_{(1\ 3)(2\ 4)}$$
$$e_{(1)}+e_{(1\ 4)}+e_{(2\ 3)}+e_{(1\ 4)(2\ 3)}$$
which, multiplied by $\beta_2$ give zero. Hence we get the equations
$$F(x_1,x_2,x_3,x_4)+F(x_2,x_1,x_3,x_4)+F(x_1,x_2,x_4,x_3)+F(x_2,x_1,x_4,x_3)=0$$
$$F(x_1,x_2,x_3,x_4)+F(x_3,x_2,x_1,x_4)+F(x_1,x_4,x_3,x_2)+F(x_3,x_4,x_1,x_2)=0$$
$$F(x_1,x_2,x_3,x_4)+F(x_1,x_3,x_2,x_4)+F(x_4,x_2,x_3,x_1)+F(x_4,x_3,x_2,x_1)=0$$
which, together with the equation produced by the second hyperplane, namely
$$\sum_{\sigma\in S_4}sgn(\sigma)F(x_{\sigma(1)},x_{\sigma(2)},x_{\sigma(3)},x_{\sigma(4)})=0$$
characterize the elements in this type of symmetry, coinciding with the equations of case 2 of page 851 of \MRS.

3) The case of functions of type $F=F_3$ is completely analogous to the previous one. Although we could take the hyperplane $x-2=0$, we will use instead the hyperplanes $-x+w+3=0$ and $x+y+w+z+1=0$. Reasoning as above, we get the characterization of this type of symmetry given by the equations
$$F(x_1,x_2,x_3,x_4)-F(x_2,x_1,x_3,x_4)-F(x_1,x_2,x_4,x_3)+F(x_2,x_1,x_4,x_3)=0$$
$$F(x_1,x_2,x_3,x_4)-F(x_3,x_2,x_1,x_4)-F(x_1,x_4,x_3,x_2)+F(x_3,x_4,x_1,x_2)=0$$
$$F(x_1,x_2,x_3,x_4)-F(x_1,x_3,x_2,x_4)-F(x_4,x_2,x_3,x_1)+F(x_4,x_3,x_2,x_1)=0$$
$$\sum_{\sigma\in S_4}F(x_{\sigma(1)},x_{\sigma(2)},x_{\sigma(3)},x_{\sigma(4)})=0$$
which are precisely the equations of case 4 of pages 852-853 of \MRS.

\bigskip

We are now in conditions of giving an explicit decomposition into symmetry classes:

\proclaim Theorem \decompositionofforms. There is a decomposition $C^d(X,K)=\bigoplus_{V\in Irr(S_d)}W^d_V(X,K)$, where $W^d_V(X,K):=\beta_VC^d(X,K)$ or, equivalently, $W^d_V(X,K)$ is the subspace of functions $F\in C^d(X,K)$ satisfying one of the following conditions:
\item{(i)} For each partition $\lambda\in P_d$ different from $(1,\dots,1)$ (i.e. the partition corresponding to the identity), we have
$$\sum_{\sigma\in\lambda}F(x_{\sigma(1)},\dots,x_{\sigma(d)})=a_{V,\lambda}F(x_1,\dots,x_d)$$
where $a_{V,\lambda}$ is the $\lambda$ coordinate of the point $a_V$.
\item{(ii)} $\beta_VF=F$.

\noindent{\it Proof.} This is nothing but Proposition \basisevaluation, in which we take $W=C^d(X,K)$. Observe that condition (ii) is equivalent to $(1-\beta_V)F=0$, and $1-\beta_V$ is the sum of all the other centrally primitive idempotents $\beta_{V'}$.
\qed

\noindent{\bf Example \examplessymmetries.} Coming back to the examples of Propositions \idempotentssymmetric\ and \coordinatessymmetric, we have:

1) If $V_s$ is the trivial representation of $S_d$, then $W_{V_s}^d(X,K)$ is the space of symmetric functions.

2) If $V_a$ is the alternating representation of $S_d$, then $W_{V_a}^d(X,K)$ is the space of skew-symmetric functions.

3) If $V_J$ is the standard representation, then $W_{V_J}^d(X,K)$ is the space of  functions with a particular type of symmetry. When $d=3$ is the $J$-symmetry of Example \Jsymmetry, when $d=4$ is the symmetry of type 3) in Example \quadrilinearrepresentations, and, in general, this is what is called the Schur functor of type $(d-1,1)$.

We apply now the techniques of this section to different examples, to show how to use the equations of the different symmetry classes that we obtained. 

\proclaim Lemma \orthogonalsymmetric. With the notation of Example \examplessymmetries, for any $F\in\bigoplus_{V\in Irr(S_d)\setminus\{V_s\}}W^d_V(X,K)$ and any $x\in X$, then $F(x,\dots,x)=0$.

\noindent{\it Proof.} The subspace $\bigoplus_{V\in Irr(S_d)\setminus\{V_s\}}W^d_V(X,K)$ is the subspace of $C^d(X,K)$ that is killed by $\beta_{V_s}=\sum_{\sigma\in S_d}e_\sigma$. Hence the elements $F\in C^d(X,K)$ in $\bigoplus_{V\in Irr(S_d)\setminus\{V_s\}}W^d_V(X,K)$ are characterized by $\sum_{\sigma\in S_d}F(x_{\sigma(1)},\dots,x_{\sigma(d)})=0$ for any $x_1,\dots,x_d\in X$. Taking $x_1=\dots=x_d=x$, we get the result.
\qed

As a proof that our reduced equations of the symmetries are simpler than what was known in the literature, we recover now the following result of Alicia Tocino (see \Tocino), who used the equations given in \MRS\ for each of the different $W_V^d(X,K)$:

\proclaim Proposition \hyperdeterminant. With the notation of Example \examplessymmetries, for any $F\in\bigoplus_{V\in Irr(S_d)\setminus\{V_s,V_J\}}W_V^d(X,K)$ and any $x,y\in X$, then 
$$F(y,x,\dots,x)=\dots=F(x,\dots,x,y)=0.$$ 
In particular, when $X$ is a vector space and $C^d(X,K)$ is the set ${X^*}^{\otimes d}$ of multilinear forms, the hyperdeterminant of any multilinear form $F$ in $\bigoplus_{V\in Irr(S_d)\setminus\{V_s,V_J\}}W_V^d(X,K)$ is zero.

\noindent{\it Proof.} A function $F\in C^d(X,K)$, is in $\bigoplus_{V\in Irr(S_d)\setminus\{V_s,V_J\}}W_V^d(X,K)$ if and only if
$$\sum_{\sigma\in S_d}F(x_{\sigma(1)},\dots,x_{\sigma(d)})=0$$
$$\sum_{\sigma\in S_d}|{\rm Fix}(\sigma)|F(x_{\sigma(1)},\dots,x_{\sigma(d)})=0$$
for any $x_1,\dots,x_d\in X$. In order to prove the proposition, we will show for example $F(x,\dots,x,y)=0$, the remaining cases being symmetric. We apply the above equalities for $x_1=\dots=x_{d-1}=x$ and $x_d=y$, dividing the summands depending on whether $\sigma(d)=d$ or not and setting $A:=F(x,\dots,x,y)$, $B:=F(y,x,\dots,x)+\dots+F(x,\dots,x,y,x)$, and we get
$$(d-1)!A+(d-1)!B=0$$
$$\big(\sum_{\sigma(d)=d}|{\rm Fix}(\sigma)|\big)A+\big(\sum_{\sigma(d)=d-1}|{\rm Fix}(\sigma)|\big)B=0$$
(we are using that $\sum_{\sigma(d)=i}|{\rm Fix}(\sigma)|$ does not depend on $i\ne d$). So it is enough to prove 
$$\sum_{\sigma(d)=d}|{\rm Fix}(\sigma)|\ne\sum_{\sigma(d)=d-1}|{\rm Fix}(\sigma)|.$$
For that purpose, we define the bijection 
$$\{\sigma\in S_d\ |Ê\ \sigma(d)=d\}\to\{\sigma\in S_d\ |Ê\ \sigma(d)=d-1\}$$ 
by assigning to any $\sigma$ with $\sigma(d)=d$ the permutation $\sigma'$ such that 
$$\sigma'(i)=\left\{\matrix{d-1&\hbox{if }i=d\hfill\cr d\hfill&\hbox{if } i=\sigma^{-1}(d-1)\hfill\cr \sigma(i)\hfill&\hbox{if }i\ne d,\sigma^{-1}(d-1)\hfill}\right.$$
Since clearly $|{\rm Fix}(\sigma)|>|{\rm Fix}(\sigma')|$, it follow that $\sum_{\sigma(d)=d}|{\rm Fix}(\sigma)|>\sum_{\sigma(d)=d-1}|{\rm Fix}(\sigma)|$, which proves the proposition.
\qed

\noindent{\bf Remark \badidea.} It seems natural to study extensions of the above result for more repetitions of elements. For example, for each $i$, one would like to characterize for which subspaces $W_V^d(X,K)$ one has $F(x,\dots,x,x_{i+1},\dots,x_d)=0$ for all $x,x_{i+1},\dots,x_d\in X$. It is interesting to observe that, if $d=4$, except for the alternating case no other symmetry subspace of $C^4(X,K)$ (of those described in \quadrilinearrepresentations) satisfies that $F(x,x,y,y)=0$ for all of the functions of that subspace.

\bigskip
\bigskip

\bigskip
\noindent{\bf References.}
\bigskip

\item{\Burnside} W. Burnside, {\it On the Condition of Reducibility of any Group of Linear Substitutions}, Proc. London Math. Soc. {\bf 3} (1905), 430-434.

\item{\Fulton} W. Fulton, {\it Algebraic Curves. An Introduction to Algebraic Geometry},  Reprint of 1969 original, Addison-Wesley, 1989. 

\item{\FultonHarris} W. Fulton, J. Harris, {\it Representation Theory}, Springer-Verlag, 1991.

\item{\Harris} J. Harris, {\it Algebraic Geometric: A First Course}, Springer-Verlag, 1992.

\item{\Lam} T.Y. Lam, {\it A Theorem of Burnside on Matrix Rings}, Amer. Math. Monthly {\bf 105} (1998), 651-653.

\item{\LuxPahlings} K. Lux, H. Pahlings, {\it Representation Theory. A Computational Approach}, Cambridge University Press, 2010.

\item{\MRS} N. Metropolis, G.-C. Rota, J.A. Stein, {\it Symmetry Classes of Functions}. J. of Algebra {\bf 171} (1995), 845-866.

\item{\Quintanilla} P. Quintanilla, {\it Representaci\'on de Grupos Finitos}, Bachelor Thesis, Facultad de Ciencias Matem\'aticas, Universidad Complutense de Madrid, July 2018 (available online at http://blogs.mat.ucm.es/arrondo/wp-content/uploads/sites/75/2023/02/TFG-Patri.pdf).

\item{\Tocino} A. Tocino, {\it The hyperdeterminant vanishes on all but two Schur functors}, J. of Algebra {\bf 450} (2016), 316-322.

\bigskip
\centerline{Instituto de Matem\'atica Interdisciplinar}
\centerline{Departamento de \'Algebra, Geometr\'{\i}a y Topolog\'{\i}a}
\centerline{Facultad de Ciencias Matem\'aticas}
\centerline{Universidad Complutense de Madrid}
\centerline{28040 Madrid, Spain}
\centerline{arrondo@mat.ucm.es}

\end